\pdfoutput=1
\documentclass[12pt]{article}
\usepackage[frozencache,cachedir=.]{minted}
\usepackage{bbm}
\usepackage{graphicx}
\usepackage{hyperref}
\hypersetup{
    colorlinks=true,
    linkcolor=blue,
    }
\usepackage{url}
\usepackage{caption}
\usepackage[dvipsnames]{xcolor}
\usepackage{enumerate}
\usepackage{amsmath,amsthm,amssymb}
\usepackage{fullpage}
\usepackage{titlesec}
\usepackage{comment}
\usepackage{mathtools}
\newcommand{\gf}[1]{\textbf{\textcolor{blue}{Gary: #1}}}
\newcommand{\gfadd}[1]{{\textcolor{red}{#1}}}
\newcommand{\aleks}[1]{\textbf{\textcolor{Orange}{Aleks: #1}}}
\newcommand{\new}[1]{\textcolor{black}{#1}}

\theoremstyle{definition}

\newtheorem{theorem}{Theorem}
\newtheorem{algorithm}{Algorithm}

\newcommand{\Prob}{{\mathrm{Prob}}}

\title{Identifying the onset and decay of \\ quasi-stationary families of almost-invariant sets  \\ with an application to atmospheric blocking events}
\date{\today}

\author{Aleksandar Badza and Gary Froyland \\ School of Mathematics and Statistics \\ UNSW Sydney NSW 2052, Australia}

\begin{document}

\maketitle

\begin{abstract}
Macroscopic features of dynamical systems such as almost-invariant   sets and coherent sets provide crucial high-level information on how the dynamics organises phase space. 
% Almost-invariant and coherent sets share the common property that  almost-invariant sets remain fixed in space, coherent sets  Classical techniques 
We introduce a method to identify \textit{time-parameterised families of almost-invariant sets} in time-dependent dynamical systems, as well as the families' \textit{emergence and disappearance}.
In contrast to coherent sets, which may freely move about in phase space over time, our technique focuses on families of metastable sets that are \textit{quasi-stationary in space}.
Our straightforward approach extends successful transfer operator methods for almost-invariant sets to time-dependent dynamics and
%identifying \emph{finite-time quasi-stationary families of almost-invariant sets}, including their appearance and breakup.
utilises the Ulam scheme for the generator of the transfer operator on a time-expanded domain. 
The new methodology is illustrated with an idealised fluid flow and with atmospheric velocity data.
We identify atmospheric blocking events in the 2003 European heatwave and  compare our technique to existing geophysical methods of blocking diagnosis.

\paragraph{Significance statement:} 
Diagnosing the emergence of persistent structures in turbulent nonlinear flows is important for many scientific and engineering applications. Computational and theoretical challenges include detecting the location of almost-invariant structures at specific times, piecing together the diagnosed objects across time while simultaneously verifying they are approximately fixed in space, and determining when these structures appear and cease to exist. Our novel spacetime approach based on an inflated generator overcomes these challenges in a single eigencomputation, and relies only on the underlying time-dependent velocity field. We also outline a simpler averaged generator to provide a first-cut diagnosis of persistence in space.
\end{abstract}

%\newpage
%\tableofcontents

%\newpage

%\section{Publication Goals and Targets}

%In this section, we detail a comprehensive plan and list targets for the publication of our findings within the coming months.

\section{Introduction}

Lagrangian approaches to analysing complex flow systems have had great success over the last two decades. Very broadly speaking these approaches may be categorised as (i) those that linearise the nonlinear dynamics directly \cite{RKLW90,PY93,HP98,Haller00,Shaddenetal05,HBV13,MB14} in phase space, through the study of spatial derivative matrices or the Cauchy--Green tensor, and (ii) those that linearise the nonlinear dynamics by moving to a space of scalar-valued functions defined on the phase space, and studying the transfer operator (or cocycles thereof) \cite{DJ99,SchuetteHabilitation,Rowley09,FSM10,FLS10,MB14,FRS19}.
The present work takes the second route and extends one of the earliest constructions in this approach, namely almost-invariant sets \cite{DJ99, SchuetteHabilitation,F05,FP09,FJK13}.
As the name suggests, almost-invariant sets are those sets that remain approximately fixed in space under the dynamics.
They have been an important forerunner to the concept of coherent sets \cite{FSM10,F13}, which share the approximate trajectory-trapping property of almost-invariant sets, but are able to rapidly move about in phase space as the dynamics evolve. 
Much progress has been made on coherent set analysis over the last decade and one of the aims of this work is to update the concept of almost-invariant sets by incorporating some of these advances.
This aim is not purely academic. 
There are real-world situations, such as the atmospheric blocking events we study here, where one is interested in identifying time-parameterised families of metastable regions \textit{that remain approximately stationary in space over a significant time duration}.
This latter criterion is a restriction that is not present in the definition of coherent sets. 

Let $M\subset \mathbb{R}^d$ be a smooth full-dimensional compact submanifold of $\mathbb{R}^d$ and let $m$ represent $d$-dimensional volume on $M$.
We assume that the dynamics on $M$ is governed by a smooth time-dependent vector field $v:[0,\tau]\times M\to\mathbb{R}^d$, i.e.\ that trajectories satisfy $\dot{x}(t)=v(t,x)$ for $t\in [0,\tau], x\in M$.
At each $t\in[0,\tau]$ we seek one or more almost-invariant sets $A_t\subset M$, where the notion of almost-invariance is made precise in the following section.
Collectively, across time, these individual almost-invariant sets form one or more time-parameterised families $\{A_t\}_{t\in[0,\tau]}$, which we require to be \textit{quasi-stationary} in the sense that the elements of this family do not vary much across time;  this will be made more precise in Section \ref{sec:qsfamilies}.
We introduce a theoretical and numerical methodology for identifying quasi-stationary families of almost-invariant sets within a long time interval $[0,\tau]$.
Of particular interest is that the family may only exist for times in a strict subset of $[0,\tau]$.
That is, our technique also identifies when these quasi-stationary families appear and disappear.

Our approach is based on generators of the Perron--Frobenius operators for flows governed by the vector fields $v(t,x)$ for various fixed $t\in[0,\tau]$.
For steady flows, the leading eigenfunctions of single generators effectively identify almost-invariant sets \cite{FJK13}.
In our time-dependent setting, we connect different generators across time using ideas from recent work on dynamic Laplacians on time-expanded domains \cite{FK23, AFK24}.
This leads to a single \textit{inflated generator} acting on spacetime, whose leading eigenfunctions encode the onset and decay of quasi-stationary families of almost-invariant sets.
Our main contributions include creating and formalising the notion of quasi-stationary families of almost-invariant sets, demonstrating quantitative relationships between the eigenspectrum of the inflated generator and almost-invariance properties of quasi-stationary families, and designing an efficient computational approach.
A further contribution is our provision of a new approach to diagnosing atmospheric blocking events, which solely relies on the use of wind fields rather than other geophysical proxies.
A blocking event, or block for short, is a quasi-stationary, anti-cyclonic %\gf{why is it anti-cyclonic? 
% is that just in the northern hemisphere?} %\aleks{Blocking events have been explicitly referred to as anti-cyclones in a lot of the literature, see for example \cite{Rex50I,Thomas21}, and tend to occur as a result of the build up of quasi-stationary high pressure systems. Anti-cyclones circulate wind around regions of high atmospheric pressure, while cyclones circulate wind around regions of low pressure. Anticyclones rotate clockwise in the Northern Hemisphere and anti-clockwise in the Southern Hemisphere. If a cyclone was to exhibit similar behaviour to a block, it is referred to as a ``cut-off low" (see https://learningweather.psu.edu/node/79). A ``cut-off low" can occur in conjuction with a block (http://www.bom.gov.au/climate/about/australian-climate-influences.shtml?bookmark=cutofflow), but to my knowledge these are seldom (if ever) referred to as blocks, or ``blocking lows".} 
air mass formed from the splitting of westerlies at mid-latitude levels in the Northern Hemisphere. 
In the Southern Hemisphere, blocks are analogously defined as cyclonic rather than anti-cyclonic air masses.
%\gf{are westerlies split or easterlies?} \aleks{Westerlies are split in the Southern hemisphere as well as the Northern hemisphere. After the westerlies are split, the branches formed drift towards the North (South) pole in the Northern (Southern) hemisphere.}. 
This trapping of air parcels inhibits regular atmospheric flow and can result in increased and prolonged extreme heat \cite{Dole83,Barriopedro10,DunnSigouin13,Thomas21}, extreme cold \cite{Dole83,Buehler11,Park11,Cai20}, or reduced air quality through increased concentration of pollutants \cite{Barnes12,Yun19,Cai20}. 
%\gf{Aleks to ensure that all citation lists are in forward chronological order.} 
Despite this phenomenon being conceptualised decades ago \cite{Rex50I,Rex50II}, much is still not known about the underlying physics behind blocks.
A better understanding of blocking physics is crucial for developing clearer links between extreme weather and blocks, and will lead to better forecasting of and preparedness for these extreme weather events \cite{Liu94,Steinfeld19,Ehstand21}.
%While blocks are not necessarily responsible for all extreme weather events, 
%Developing such an understanding into blocks begins by gaining an idea of how these blocks take shape over the course of their lifespan, specifically how a block forms, .

Prior research into blocks in the geophysics literature has often focused on statistics;  that is, blocking frequencies at particular lines of latitude.
%\cite{Lejenas83,Tibaldi90,Sausen95,Pelly03,Scherrer06,Barriopedro10,DunnSigouin13,Lucarini20,Lupo21,Wazneh21} . 
These frequencies are estimated using blocking proxies, such as geopotential height gradients \cite{Lejenas83,Tibaldi90,Tibaldi95,Scherrer06,Pinheiro19,Thomas21}, geopotential height anomalies \cite{Dole83,Sausen95,Barriopedro10,Pinheiro19,Thomas21} or potential vorticity anomalies \cite{Schwierz04,Scherrer06,Pinheiro19,Thomas21}. 
This frequency analysis largely prohibits further examination of the key blocking physics and impacts on weather patterns. 
More recently, post-analysis methods investigating the dynamics of blocks that have already been determined from geophysical proxies have appeared \cite{Hirt18,Steinfeld19,Hauser23,Schoeller24}.
%Dynamical systems methods of block identification that use one or more observed geophysical inputs include \cite{Faranda15, SerGiacomi15, Schubert16, Nakamura18, Ehstand21}.

A range of dynamical systems methods for block detection have also been developed and implemented on observational data  \cite{Faranda15,SerGiacomi15,Nakamura18,Ehstand21,VonLindheim21} or simulated data obtained from barotropic vorticity models \cite{Shutts83,Tantet15,Schubert16,Ghil18,Xavier24}.
\new{In the latter case blocks are typically identified as persistent, large-scale anomalies of tropospheric air circulation (a process driven by a large number of atmospheric and thermodynamic variables \cite{Ghil18}) which last at least a few days but do not span the entire length of a meteorological season.}  
The methods above either (i) capture the two-dimensional extent of %\gf{illustrate means ``capture the two-dimensional extent''?}
blocking events at a single instance of time (in particular without capturing the full evolution of a block over its lifespan)
%\gf{if these methods capture the two-dimensional extents at different times, why can't they be put together to ``capture the full evolution''?}
\cite{Shutts83,SerGiacomi15,Schubert16,Ghil18,Ehstand21,Xavier24}, (ii) identify longitude levels (but not two-dimensional regions) which remain blocked over a period of time \cite{Nakamura18}, (iii) detect the timing of the transition of an atmospheric flow regime to a blocked state and don't provide information on the spatial location of the block \cite{Faranda15,Tantet15,Ghil18}, or (iv) ignore the quasi-stationarity property of blocks \cite{Shutts83,SerGiacomi15, Ehstand21,VonLindheim21,Xavier24}.

The aim of our new approach is to provide a full description of a blocking event over its lifetime, by pinpointing the block's onset and decay in time, and tracking its spatial extent throughout its lifetime.
%accomplish precisely this, without any reliance on the aforementioned geophysical block indication proxies or other geophysical quantities outside of wind fields.
%The aim of our new approach is to produce a comprehensive visualisation of a block from its onset right up until its decay, without any reliance on the aforementioned geophysical block indication proxies or other geophysical quantities outside of wind fields.
Our technique directly uses wind velocities and does not require other geophysical proxies.
We explicitly incorporate the crucial quasi-stationarity property of blocks \cite{Rex50I,Rex50II,Shutts83,Hirt18,Steinfeld19,Kautz22,Hauser23}, a property which is often of peripheral importance in frequency analysis \cite{Rex50I,Lejenas83,Tibaldi90,Scherrer06,Barriopedro10,DunnSigouin13} but increasingly used in other dynamical systems based approaches \cite{Faranda15,Tantet15,Schubert16,Nakamura18,VonLindheim21}.
To demonstrate this new approach, we study two blocking events from the Northern Hemisphere summer of 2003, during what is regarded as some of the worst extreme heat events recorded in Europe \cite{Black04,Ehstand21,Thomas21}.

The outline of this paper is as follows. In Section \ref{sec:2}, we recap relevant theory concerning generators of steady flows from \cite{FJK13}.  This section ends with a bridging construction of an \textit{averaged generator}.
In Section \ref{sec:qsfamilies} we construct the inflated generator, inspired by the constructions in \cite{FJK13,FK23}, and detail how this operator is used to identify quasi-stationary families of almost-invariant sets. In Section \ref{sec:disc}, we describe the numerical implementation of the inflated generator method, including the discretisation of the generator, selection of key parameters for the generator, and the extraction of time-parameterised families of almost-invariant sets from its spacetime eigenfunctions. 
In Section \ref{IdealisedDGFlow}, we illustrate our theory on an idealised fluid flow that possesses well-defined quasi-stationary almost-invariant sets. 
In Section \ref{sec:block}, we apply the inflated generator method to European Centre for Medium-Range Weather Forecasts (ECMWF) atmospheric velocity data to identify two key blocking events that occurred in the European summer of 2003 \cite{Black04,Thomas21}.
We demonstrate that our inflated generator method produces a clear, detailed visualisation of the evolution of a block from its onset right up until its decay. 
We show that identical results are not obtained from three frequently-used geophysical block detection methods, partly because of their reliance on blocking proxies defined through scalar fields on individual time slices, without explicit consideration of quasi-stationarity of the block.
%, which are only loosely connected between each time slice. 
\begin{comment}
\gf{Aleks to add what we found (very briefly) in this comparison.}

\gf{This sounds pretty bland, like my contribution criticism above. Again we should focus on novel contributions}.
\end{comment}

%Coherent set techniques ought to find atmospheric blocks, but there may be many confounding false positives, namely coherent atmospheric regions that move rapidly in space and 

\section{Almost-invariant sets and generators of steady flows}
\label{sec:2}
As a first step toward our main constructions in the following section, we suppose for the moment that we have a steady flow given by the vector field $v(t,\cdot)$ at some \textit{fixed time $t$}.
To simplify the mathematical theory of spectra of our generators, we will consider the \textit{stochastic} flow given by 
\begin{equation}
    \label{sde}
dx_s=v(t,x)\, ds + \epsilon\, dw_s,
\end{equation}
where $\epsilon$ is a small positive number representing the magnitude of Brownian motion $\{w_s\}$. The corresponding forward propagator of a function $f:M\to\mathbb{R}$ at time $t$ for duration $s$ under the fixed drift $v(t,\cdot)$ is given by the Perron--Frobenius operator $P_t^s.$

Following \cite{FJK13}, given a probability measure $\nu$ on $M$ we say that a set $A\subset M$ is almost-invariant with respect to $\nu$ at time $t$ over a short time duration $s>0$ if $$\rho_{t,\nu}^s(A):=\frac{\Prob_\nu(x_t\in A, x_{t+s}\in A)}{\Prob_\nu(x_t\in A)}\approx 1,$$
where $\Prob_\nu$ refers to the distribution of points in space at time $t$ according to $\nu$.
If we wish to find almost-invariant sets for this steady flow in $M$, following the approach in \cite{FJK13} we construct the \textit{generator} of the Perron--Frobenius operator $P_t^s$ for this steady flow at this fixed time $t$:
%The generator for the corresponding Perron--Frobenius operator is:
\begin{equation}
    \label{gen}
    G_tf(x):=-\nabla\cdot(v(t,x)f(x))+\frac{\epsilon^2}{2}\Delta f(x).
\end{equation}
% \begin{equation}
%     \label{genadj}
% G_t^*f(x):=v(t,x)\cdot\nabla f(x)+\frac{\epsilon^2}{2}\Delta f(x).
% \end{equation}
To connect the generator and $\rho^s_{t,\nu}$, we define the functional 
\begin{equation}
    \label{Gtfunctional}
G_{t,A}f:=\lim_{s\to 0} \int_A \frac{P_t^sf-f}{s}\ dm,
\end{equation}
which is well defined for $f\in D_A(G_t)$ (denoting the linear subspace of $L^1$ where this limit exists).
For example, denoting the indicator on $A$ by $\mathbf{1}_A$,
Proposition 3.2  \cite{FJK13}  states
\begin{equation}
    \label{rhogen}
    \rho^s_{t,m}(A)=1+\frac{G_{t,A}\mathbf{1}_A}{m(A)}\cdot s+o(s),
\end{equation}
with $G_{t,A}\mathbf{1}_A<0$.
Thus, almost-invariant sets at time $t$ are those sets $A$ for which $G_{t,A}\mathbf{1}_A/m(A)\lessapprox 0$.
Our strategy is to use the leading eigenfunctions of $G_t$ to find almost-invariant sets at time $t$. 
We now briefly motivate this approach recapping another result from \cite{FJK13}.

First, applying zero Neumann (reflecting) boundary conditions at the boundary of $M$, we have a well-posed eigenvalue problem:
\begin{equation}
    \label{steadyeigenprob}
G_tf=\lambda f,\qquad \nabla f(x)\cdot n(x)=0\mbox{ for }x\in\partial M,
\end{equation}
where $n(x)$ is a vector normal to $\partial M$ at $x\in\partial M$.
%It is clear that the constant function $f(x)\equiv \mathbf{1}(x)$ is in the nullspace of $G_t^*$.
By standard theory, e.g.\ Theorem 8.4.4 \cite{mcowen}, the spectrum\footnote{\new{considering $G_t$ as an operator from $H^2(M)\cap H_0^1(M)$ to $L^2(M)$, where $H^2(M)$ is the Sobolev space of functions with weak second derivatives and $H_0^1(M)$ is the Sobolev space of functions with weak first derivatives, which vanish on $\partial M$.}} of $G_t$   consists of isolated eigenvalues $\{\lambda_k\}_{k=1}^\infty$ of finite multiplicity and no finite accumulation points.
We are primarily concerned with the eigenvalues of largest (negative) real part, which corresponds to the eigenfunctions with the slowest rates of decay under the action of the generator $G_t$.
%The leading eigenvalue $\lambda_1=0$ is simple.
\begin{theorem}[Theorem 3.3 \cite{FJK13}]
\label{lamthm}
Let $\lambda<0$ be a real eigenvalue for $G_t$ with eigenfunction $f$, normalised so that $\int_M |f| \ dm=1$. Then $\int_M f\ dm=0$. Define the signed measure $\nu(A)=\int_A f\ dm$, subsets $A^+=\{f\ge 0\}$, $A^-=\{f<0\}$, and suppose $f\mathbf{1}_{A^\pm}\in D_A(G_t)$.
Then 
\begin{equation}
    \label{DJ99gen}
    G_{t,A^+}(f\mathbf{1}_{A^+})/\nu(A^+)+G_{t,A^-}(f\mathbf{1}_{A^-})/\nu(A^-)=\lambda.
\end{equation}
\end{theorem}
If we choose $\lambda<0$ to be a real eigenvalue close to 0 then the sets $A^\pm\subset M$ will be almost-invariant with respect to the probability measure with density $|f|$.
Thus, one strategy for finding almost-invariant sets is to look at positive and negative parts of dominant eigenfunctions of $G_t$.
This idea goes back to \cite{DJ99} and has been refined over the years, by way of optimising thresholds to consider $A^+=\{f\ge c\}$, $A^-=\{f<c\}$ \cite{F05}, considering multiple eigenvectors \cite{FD03,F05}, and the subsequent extraction of multiple almost-invariant sets by various means, e.g.\ $c$-means clustering \cite{F05}, constrained nonlinear optimisation  \cite{Deuflhard05}, and sparse optimisation \cite{FRS19}.
We will use the latter approach of sparse eigenbasis approximation (SEBA) \cite{FRS19}, which is robust and well suited to extracting many features, particularly when the features do not fully partition the phase space.

\paragraph{An averaged generator for unsteady flows.}
We have just discussed connections between generators of steady vector fields $v(t,\cdot)$ for fixed $t$ and almost invariant sets for the steady flows.
As a partial bridge to the next section, we now consider the full time-dependent vector field $v(t,x)$ for $0\le t\le \tau$, and ask for sets in $M$ that are almost-invariant \textit{in a time-averaged sense} across $[0,\tau]$.
That is, we seek $A\subset M$ for which the \textit{average} residence fraction $$\bar{\rho}^s_m(A):=\new{\frac{1}{\tau}}\int_0^\tau \rho^s_{t,m}(A)\ dt$$ is close to 1.
Define $\bar{v}(x)=\frac{1}{\tau}\int_0^\tau v(t,x)\ dt$ to be the averaged vector field at a fixed location $x$ across the time interval $[0,\tau]$.
We now introduce the time-averaged generator:
\begin{equation}
    \label{bargen}
     \bar{G}f(x):=\frac{1}{\tau}\int_0^\tau G_tf(x)\ dt=-\nabla\cdot(\bar{v}(x)f(x))+\frac{\epsilon^2}{2}\Delta f(x).
\end{equation}
From \eqref{rhogen} we immediately obtain 
\begin{equation}
    \label{barrhogen}
    \bar\rho^s_{m}(A)=1+\frac{\bar G_{A}\mathbf{1}_A}{m(A)}\cdot s+o(s),
\end{equation}
with 
\begin{equation}
    \label{Gbarfunctional}
\bar G_{A}f:=\lim_{s\to 0} \int_A \frac{\bar P^sf-f}{s}\ dm,
\end{equation}
and $\bar P^s$ denoting the propagator (Perron--Frobenius operator) of $f$ from time $0$ to $s$ under the SDE \eqref{sde} with $v(t,x)$ replaced by $\bar v$.
An analogous version of Theorem \ref{lamthm} holds with $G_t, G_{t,A^+}, G_{t,A^-}$ replaced by $\bar G, \bar G_{A^+}, \bar G_{A^-}$, respectively.
In summary, to find subsets that are \textit{persistently almost invariant across the full finite time interval $[0,\tau]$} one can use eigenfunctions of $\bar G$ in the algorithms that follow.
Our ultimate aim is more ambitious, namely finding \textit{quasi-stationary families} of almost-invariant sets, including their \textit{emergence and disappearance}.
The averaged generator may be viewed as a simpler, less computationally intensive, first-cut analysis tool.

%The dynamic generator $\mathcal{G}^D$ acting on functions on $M$ is defined as $\mathcal{G}^D=\frac{1}{\tau}\int_0^\tau G_t\ dt$.
%\textbf{Note however that $\mathcal{G}^D$ is equal to the generator corresponding to the time-averaged vector fields...}

\section{Quasi-stationary families of almost-invariant sets and an inflated generator}
\label{sec:qsfamilies}

%\subsection{The inflated generator}

\begin{figure}[htb]
      \centering
\includegraphics[width=0.6\textwidth]{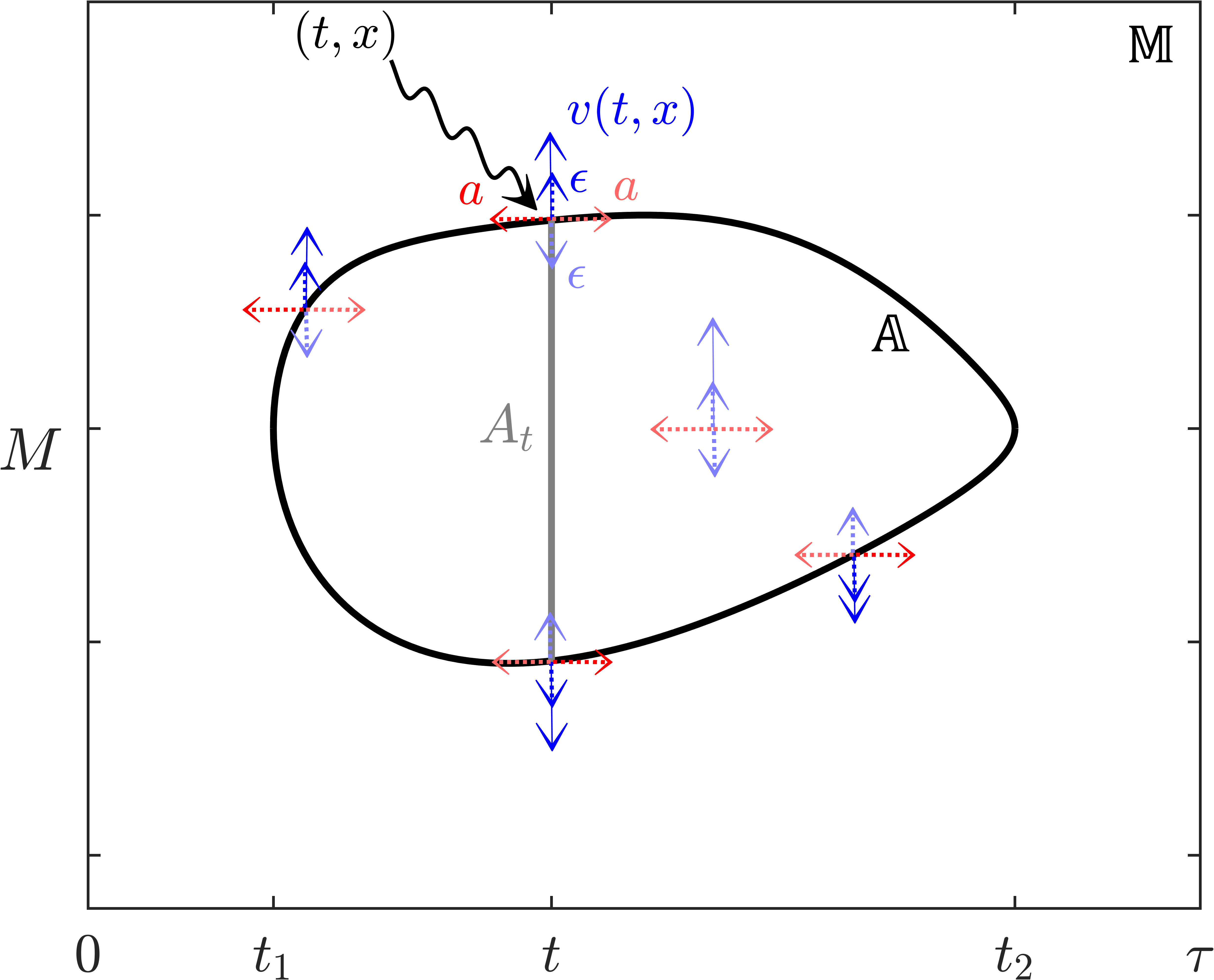}
      \caption{
      %\gf{looks better.  The $A_t$ seems to be grey instead of black. so too the $\epsilon$. It is probably good to vary the length of the drift arrows. We should also probably include $\epsilon$ arrows over the top of the $v$ arrows to make it clear that the diffusion is symmetric;  these diffusion arrows can always be shorter than the drift arrows. the $v(t,\cdot)$ is good, but it should be $v(t,x)$, and we need to indicate that the point is $(t,x)$, perhaps by pointing to it with a wiggly black arrow.} \aleks{Thank you for that, a new schematic is coming in now...}
        A schematic of a space-time set $\mathbb{A} \subset \mathbb{M}=[0,\tau]\times M$ for one-dimensional $M$.
      We seek $\mathbb{A}$ that are almost-invariant under the spacetime SDE (\ref{sdeaug}). 
    A time fibre $A_t$ of $\mathbb{A}$ is shown in grey.
      Trajectories may exit $\mathbb{A}$ through $\partial \mathbb{A}$ via drift $v$ in space and/or diffusion in space of strength $\epsilon$, indicated by the solid and dotted blue arrows respectively.
      Therefore, the boundary of $\mathbb{A}$ should attempt to trace through spacetime so that the component of $v(t,x)$ normal to $\partial\mathbb{A}$ is small, to minimise trajectory loss through drift;  \textit{this will ensure that each $A_t$ is almost-invariant}.
      Trajectories may also exit $\mathbb{A}$ through $\partial\mathbb{A}$ via diffusion in the temporal direction of strength $a$, indicated by the dotted red arrows. 
      If the time fibres vary slowly in time, the loss of trajectories due to temporal diffusion will be small;  \textit{this will ensure that the family $\{A_t\}_{t=0}^\tau$ is quasi-stationary}. Dark blue and red arrows indicate trajectory loss from $\mathbb{A}$, while pale blue and red arrows illustrate motion internal to $\mathbb{A}$.}
     % with the $t^{\rm th}$ time fibre $A_t$ highlighted in grey. In this Figure, mass loss due to drift is indicated using solid arrows, while dotted arrows represent mass loss due to diffusion. Blue arrows represent mass loss due to spatial drift/diffusion, while red arrows represent mass loss due to temporal diffusion. Arrows with paler colours represent mass movement within $\mathbb{A}$, while arrows with darker colours represent mass loss from $\mathbb{A}$.}
      \label{fig:schematic}
\end{figure}

In this section we explain how to extend the above ideas to find quasi-stationary families of almost-invariant sets in a time-dependent flow over a finite time duration $[0,\tau]$, which is one of the main contributions of this work.
We define a time-expanded domain $\mathbb{M}=[0,\tau]\times M$ and construct a stochastic process on this domain as follows.
\begin{enumerate}
    \item In the one-dimensional time coordinate we have diffusion in time with intensity $a>0$ and zero drift,
    \item In the original space coordinates, on $M$ at time $t$ (i.e.\ on $\{t\}\times M\subset \mathbb{M}$) we have diffusion with intensity $\epsilon$ and drift given by $v(t,\cdot)$.
    \end{enumerate} 
More precisely, we define points in $\mathbb{M}$ as $X=(t,x)$, a vector field $V:\mathbb{M}\to\mathbb{R}^{d+1}$ as $V(t,x)=\begin{pmatrix}
     0\\
     v(t,x)  
\end{pmatrix}$ and a $(d+1)\times (d+1)$ diffusion matrix $\Sigma=\begin{pmatrix}
    a & 0 \\ 0 & \epsilon I_d
\end{pmatrix}$, and consider the stochastic flow given by 
\begin{equation}
\label{sdeaug}
    dX_s=V(X_s)\, ds + \Sigma\, dW_s,
\end{equation}
where $dW_s$ is Brownian motion in $\mathbb{M}$.
%\aleks{Should we have $dt$ or $ds$ in the above equation?}\gf{ds. t is a coordinate, just like x. see also equation 1.}
Note that the SDE \eqref{sdeaug} is homogeneous in time due to the time extension.
The idea behind this construction is as follows.
Let $\mathbbm{m}$ denote volume on $\mathbb{M}$.
For a set $\mathbb{A}\subset\mathbb{M}$, denote the $t^{\rm th}$ time fibre $A_t=\mathbb{A}\cap(\{t\}\times M)$, which we consider as a subset of $M$;  see Figure \ref{fig:schematic}.
\begin{comment}
    \gf{the $\mathbb{A}$ and $\mathbb{M}$ are huge compared to all over notation, they should be the same size. Overall, the arrows and objects are large compared to the axis labels and other notation, so match thoese better too (could possibly be done by simply resizing the window and seeing what happens).} A schematic of a quasi-stationary family of almost--invariant sets $\mathbb{A} \subset \mathbb{M}$ \gf{$\mathbb{A}$is a spacetime set, not a family.  The family is the collection of fibres}. \gf{$A_t$ (too small) is not referenced in the caption, but is a key player in the whole figure.} In this Figure, solid blue arrows represent mass loss from $\mathbb{A}$ due to spatial diffusion \gf{where is the loss due to drift indicated in the figure?  Also, there is notation for $a$, $\epsilon$, but not for drift}, while solid red arrows represent mass loss due to temporal diffusion.\gf{I thought that drift was going to be solid and diffusion was going to be dotted.  Internal was going to be pale and external was dark.  All points should have four arrows.} Dashed arrows represent mass movement due to diffusion within $\mathbb{A}$ and have been coloured with paler tints of blue and red. The solid pale blue arrow represents spatial movement within $\mathbb{A}$ due to drift. 
\end{comment}
We may now write the space-time set $\mathbb{A}$ as the disjoint union of the sets $A_t$ placed on their individual time-$t$ fibres: 
 $\mathbb{A}=\bigcup_{t\in [0,\tau]} (\{t\}\times A_t)$.
Note that the dynamics restricted to the time-fibre $A_t$ are exactly the spatial dynamics we discussed in the previous section.
For $s>0$ consider the ratio
$$\rho_{\mathbbm{m}}^s(\mathbb{A}):=\frac{\Prob_\mathbbm{m}(X_0\in \mathbb{A}, X_{s}\in \mathbb{A})}{\Prob_\mathbbm{m}(X_s\in \mathbb{A})}.$$
If $\rho_\mathbbm{m}^s(\mathbb{A})$ is close to 1, this means that on average for a flow duration $s$:
\begin{enumerate}
    \item \textbf{Spatial almost-invariance:} the dynamics of (\ref{sdeaug}) restricted to each time fibre $A_t$ (given by \eqref{sde}) is such that trajectories rarely leave $A_t$.
    \item \textbf{Quasi-stationary family in time:} the variation of $A_t$ with $t$ is small in order that the pure diffusion in the time direction in (\ref{sdeaug}) allows only rare exits from the family $\{A_t\}_{t\in[0,\tau]}$. 
\end{enumerate}
It is important to note that the quasi-stationary family $\{A_t\}_{t=0}^\tau$ may consist of nonempty $A_t$ only for a subset of $[0,\tau]$, for example, for $t\in[t_1,t_2]$, where $0<t_1<t_2<\tau$, as illustrated in the example in  Figure \ref{fig:schematic}.
Thus by considering almost-invariant sets $\mathbb{A}$ of the SDE \eqref{sdeaug} we achieve a third property, in addition to the two above, namely we see \textit{when a quasi-stationary family of almost-invariant sets is born} (when the $A_t$ become nonempty) and \textit{when this family dissipates} (when the $A_t$ become empty again).
%These three points are illustrated in Figure \ref{fig:schematic}.

We now define a generator on $\mathbb{M}$ in order to identify such quasi-stationary families of almost-invariant sets from its eigenfunctions.
We define the \textit{inflated generator} of the Perron--Frobenius operator $\mathcal{P}$ for the flow generated by the SDE \eqref{sdeaug} by
% \begin{equation}
%     \label{infgen}
% \mathcal{G}F(t,x)=G_tF(t,x)+a^2\partial_t^2F(t,x),
% \end{equation}
% with adjoint
% \begin{equation}
%     \label{infgenadj}
% \mathcal{G}^*F(t,x)=G_t^*F(t,x)+a^2\partial_t^2F(t,x).
% \end{equation}

%we may write the inflated generator and its adjoint as generators acting on functions in spacetime.
\begin{equation}
    \label{infgenasgen}
    \mathcal{G}F(t,x)=-\nabla\cdot(V(t,x)F(t,x))+\frac{\Sigma^2}{2}\Delta F(t,x),
\end{equation}
where $F:\mathbb{M}\to\mathbb{C}$.
In analogy with \eqref{Gtfunctional}, define the functional 
\begin{equation}
    \label{Gfunctional}
\mathcal{G}_\mathbb{A}F:=\lim_{s\to 0}\int_\mathbb{A} \frac{\mathcal{P}^sF-F}{s}\ d\mathbbm{m}.
\end{equation}
%\aleks{Should we have $\mathcal{P}^s_0$ above?}\gf{no, we are in spacetime so the operator is time independent. i removed the zero from the averaged operator to be consistent, and i slightly altered the description of $P_t^s$ earlier to equivalently say "duration $s$" rather than "to time $t+s$" to emphasise the duration aspect.} 
Because $\mathcal{G}$ is a generator of a homogeneous process we immediately obtain (again from Proposition 3.2 \cite{FJK13})
\begin{equation}
    \label{rhogenaug}
\rho_\mathbbm{m}^s(\mathbb{A})=1+\frac{\mathcal{G}_{\mathbb{A}}\mathbf{1}_\mathbb{A}}{\mathbbm{m}(\mathbb{A})}\cdot s+o(s).
\end{equation}
Thus, almost-invariant sets $\mathbb{A}$ satisfy $\mathcal{G}_{\mathbb{A}}\mathbf{1}_\mathbb{A}/\mathbbm{m}(\mathbb{A})\lessapprox 0$.
We apply zero Neumann boundary conditions at the boundary of $\mathbb{M}$ to obtain the eigenvalue problem:
\begin{equation}
    \label{spacetimeeigenprob}
\mathcal{G}F=\Lambda F,\qquad \nabla F(x)\cdot N(x)=0\mbox{ for }X\in\partial \mathbb{M},
\end{equation}
where $N(X)$ is a vector normal to $\partial \mathbb{M}$ at $X\in\partial \mathbb{M}$.
%It is clear that the constant function $f(x)\equiv \mathbf{1}(x)$ is in the nullspace of $G_t^*$.
As in the previous subsection, the spectrum of $\mathcal{G}$ consists of isolated eigenvalues $\{\Lambda_k\}_{k=1}^\infty$ of finite multiplicity and no finite accumulation points, and we are interested in the eigenvalues with largest (negative) real part.
It is not hard to see that an analogous version of Theorem \ref{lamthm} holds.

In practice we will work with the $L^2$-adjoint of the inflated generator:
\begin{equation}
\label{infgenadjasgen}
\mathcal{G}^*F(t,x)=V(t,x)\cdot\nabla F(t,x)+\frac{\Sigma^2}{2}\Delta F(t,x).
\end{equation}
This is consistent with the use of right eigenvectors of Markov chain approximations of Perron--Frobenius operators to identify almost-invariant sets \cite{F05,FP09}.
In \cite{F05,FP09}, the transition matrix $P$ of a Markov chain is row-stochastic, requiring vector multiplication on the left to evolve forward in time.
Thus, right eigenvectors\footnote{Strictly speaking, in \cite{F05,FP09} a symmetrised transition matrix $R$ was constructed from the non-symmetric transition matrix $P$ and the right eigenvectors of $R$ were used. Nevertheless the same ideas apply to the right eigenvectors of the transition matrix $P$ itself, without the optimality properties afforded by the symmetrisation.} of $P$ correspond to left eigenvectors of $P^*$, and therefore to eigenfunctions of $\mathcal{G}^*$ in our current context. 

A second reason to consider $\mathcal{G}^*$ is to produce robust results when flows have regions of convergence in $M$.
If the vector field $v(t,x)$ has stable approximately fixed points (or more general attracting sets) then the eigenfunctions of $G_t$ will tend to localise about these.
On the other hand, the leading eigenfunction of $G_t^*$ -- corresponding to the eigenvalue 0 -- is always $f_1=\mathbf{1}_M$, and the remaining leading eigenfunctions of $G_t^*$ tend to have a strong level set structure, where each approximate level set indicates an approximate basin of attraction of a given stable nearly-fixed point.
This is explained in more detail in section V \cite{FSvS14}, in the context of  finite-state Markov chains  with multiple absorbing regions and applied to map plastic garbage patches in the ocean. 
There the transition matrices $P$ of the Markov chains are again row-stochastic and right eigenvectors are used, which corresponds to left eigenvectors of $P^*$ (and therefore to eigenvectors of $G_t^*$ in our current context). 
See also \cite{K11} for a related approach based on classical formulae for absorption probabilities in the presence of a single absorbing state.

The spectrum of $\mathcal{G}^*$ also consists of isolated eigenvalues $\{\Lambda_k\}_{k=1}^\infty$ of finite multiplicity and no finite accumulation points.
The eigenvalue $\Lambda_1=0$ has the corresponding constant eigenfunction $F_1^*=\mathbf{1}_\mathbb{M}$.
Because $\mathcal{G}^*$ is the adjoint of the generator of a stochastic process, the same considerations in \cite{F05,FP09,FSvS14} apply, and leading eigenfunctions $F_1^*, F_2^*,\ldots$ of $\mathcal{G}^*$ will tend to have level set structures on almost-invariant sets\footnote{When attracting sets exist, we expect level set structures on basins of attraction of almost-invariant sets.} in $\mathbb{M}$.

The space-time block structure of $V$ and $\Sigma$ allows us to immediately identify some of the eigenfunctions of $\mathcal{G}^*$.
Following the arguments in Section 3.4.2 \cite{FK23}, we note that the constant-in-space functions $F_k^*(t,\cdot):=\cos(k\pi t/\tau)$ are eigenfunctions of $\mathcal{G}^*$ with eigenvalues $\Lambda_k^{\rm temp}:=-a^2\pi^2k^2/2\tau^2$.
We call the $F_k^*$ \textit{temporal eigenfunctions} because they only vary in time and carry no spatial information.
%because the only relevant component of $\mathcal{G}^*$ is the one-dimensional ``time'' Laplacian term %\textbf{note problem with $a$ vs $\epsilon$.}
We typically want to choose the temporal diffusion strength $a>0$ large enough so that there is only around one temporal eigenvalue appearing before we begin to see the non-temporal (spatial) part of the spectrum, whose eigenvalues contain useful spatiotemporal information relevant to the dynamics of $\mathbb{M}$.  

% \begin{enumerate}
%     \item Spectrum?  In the left-half plane? (see Pavliotis book)
%     \item Elliptic operators?
%     \item Spectrum is discrete?
% \end{enumerate}

Finally, we note that if $F\in L^2(\mathbb{M})$ is an eigenfunction of $\mathcal{G}$ corresponding to a non-temporal (spatial) eigenfunction of $\mathcal{G}^*$ then $F$ is orthogonal to all $F_k^*$. 
Let $W_0=\{f\in L^2([0,\tau]): \int_0^\tau f(t)\ dt=0\}=\mathbf{1}^\perp_{[0,\tau]}$.
Since the temporal eigenfunctions are dense in $\{f\mathbf{1}_M: f\in W_0\}\subset L^2(\mathbb{M})$, we may follow the argument of Section 4.2 \cite{FK23} to see that $\int_M F(\cdot,x)\ d\ell(x)=\mathrm{const}\  \mbox{a.e.\ on $[0,\tau]$}$.
This provides a numerical test for whether $F_k^*$ is a temporal eigenfunction which involves computing the variation of the spatial integrals of the corresponding $F_k$.
In practice it is often clear by inspection which eigenfunctions $F_k^*$ are temporal.

%Regarding the eigenvectors of $\mathcal{G}^*$ for use in identifying almost-invariant sets, we in practice, we will calculate $\mathcal{G}^*$ for the negative vector field $-v(t,x)$ so that the eigenvectors of the adjoint operator correspond to forward-time motion.

%\textit{It is probably better (unfortunately) to compute left eigenvectors (i.e. efuns of $\mathcal{G}$ and to divide these by the leading efun of $\mathcal{G}$.  Then we obtain time slices that are forward-time derivatives and they look more dynamically consistent when viewed as a movie in forward time. These new functions have the property that they are orthogonal to the leading efun of $\mathcal{G}$, so that they lack the constant time-slice property like the efuns of $\mathcal{G}^*$ and if one wants to identify temporal efuns they should first compute efuns of $\mathcal{G}^*$.  This is a little unfortunate because of the extra complexity.
%I have tested this on double gyre and it should be also checked on the ecmwf data.}

%Quote and adapt material from section 3.2 \cite{FJK13} about almost-invariance and its approximation with the generator for short times (Prop 3.2 \cite{FJK13}).
%In the adaptation, this would be in spacetime;  I am not sure how clean/possible this is with the time diffusion.
%Try to adapt Theorem 3.3 \cite{FJK13} as well.

%Maybe something can also be done with escape rates in section 3.5 \cite{FJK13}.  At least it could be used as heuristic motivation if a true spacetime generalisation is difficult.

\section{Discretisation of generators}
\label{sec:disc}
%In this section, we detail the computational steps undertaken to obtain results for our Dynamic Laplacian analysis on our flow systems of interest. %The section ideally should contain steps laid out with formulae and step by step instructions, though algorithms and code can also be inserted in the below sections if desired/appropriate.

This section concerns the numerical discretisation of $\mathcal{G}$.
Because of the partially separable structure of $\mathcal{G}$ in space and time evident in \eqref{infgen}, we discretise time and space separately.
The time interval $[0,\tau]$ is sampled with equispaced times $0=t_0<t_1<\cdots<t_T=\tau$ with $t_l=l\tau/T$, $l=0,\ldots,T$.
The manifold $M$ is discretised into boxes $B_1,\ldots,B_N\subset\mathbb{R}^d$ with common side lengths $\ell_1,\ldots,\ell_d$.
For each $t_l$, $l=0,\ldots,T$, 
we discretise the copy $\{t_l\}\times M$ of the manifold $M$.
This $t_l$-fibre is simply a copy of $M$ on which the vector field is $v(t_l,\cdot)$.

\subsection{Discretisation of spatial generators on each time fibre}
\label{sec:gendisc}
For simplicity, we use the ``Ulam's method for the generator'' scheme elaborated in \cite{FJK13}.
The main construction is given in \cite[Equation (4.4)]{FJK13}, which involves a simple spatial integration on the faces of the cubes discretising $M$ with respect to $d-1$-dimensional volume, denoted $m_{d-1}$.
For the computation on a single time-fibre $\{t_l\}\times M$, the integrand on a particular cube face $B_i\cap B_j$ is the dot product of the vector field $v(t_l,\cdot)$ with the outward-pointing unit normal vector $n_{ij}$, representing the rate of flux out of the face.
With $N$ cubes ($\ell_1=\ell_2=\cdots=\ell_d$), one obtains a sparse $N\times N$ matrix $\mathbf{G}^{t_l}$ where $\mathbf{G}^{t_l}_{ij}$ is the proportional rate of outflow from cube $B_i$ into $B_j$. 
\begin{equation}
    \label{Geqn}
    \mathbf{G}^{t_l}_{ij}=\begin{cases}
    \frac{1}{m(B_j)}\left(\int_{B_i\cap B_j}(\max\{v(t_l,x)\cdot n_{ij},0\})\ dm_{d-1}(x)\right) + \epsilon^2/2\ell^2,&\quad i\neq j,\\
    -\sum_{j\neq i}\mathbf{G}^{t_l}_{ij},&\quad\mbox{otherwise.}
\end{cases}
\end{equation}
The term $\epsilon^2/2\ell^2$ in (\ref{Geqn}), not present in \cite{FJK13}, is a central-difference estimate of $\epsilon^2/2$-strength diffusion that we add to ensure a simple leading eigenvalue.  
The small spatial diffusion parameter $\epsilon$ is set to be \begin{equation}
    \label{epseqn}
\epsilon=\sqrt{0.1\bar{v}\ell}
\end{equation}
where\footnote{Since $\epsilon^2/2\ell^2=0.05\bar{v}/\ell$, one may view this diffusion term as adding a drift of 5\% of the median velocity norm $\bar{v}$ to the integrand in \eqref{Geqn}.}  $\bar{v}$ is the median of $\|v(t,x)\|_2$ over $\mathbb{M}$. 
%\gf{GF to convince himself that \eqref{Geqn} is what one would obtain if one added $\epsilon$-diffusion through central differencing.}
%and is added to the outflow across each face, consistent with a second-order central-difference estimate of $(\epsilon^2/2)\Delta$.
%centralto ensure that each $\mathbf{G}^{t_i}$ has a one-dimensional kernel.
The only nonzero entries of $\mathbf{G}$ are those where $B_i$ and $B_j$ share a face.
If $\mathbf{g}\in\mathbb{R}^N$ is a vector whose entries represent the mass contained in each cube, then the matrix-vector product $\mathbf{g}\mathbf{G}$ is a vector containing the rates of change of the masses in each cube.
In our experiments we use julia's \verb"QuadGK.jl" package \cite{quadgk} to estimate the integrals in \eqref{Geqn} with Gauss-Kronrod quadrature of order 1 with absolute and relative tolerance of $10^{-2}$. 
In section \ref{EuroBlock} we generalise \eqref{Geqn} to the case of non-identical non-cubic partitions.

%This construction is applied to each time fibre $\{t_i\}\times M$ with the vector field $v(t_i,\cdot)$ to yield a matrix $\mathbf{G}^{t_i}$.

%\subsection{Discretisation of the dynamic Generator}

%The discretised dynamic generator is simply the average of the discretised generators $\mathbf{G}^D=\frac{1}{T+1}\sum_{i=0}^T \mathbf{G}^{t_i}$.
%In julia code, if the $\mathbf{G}^{t_i}$ are stored in a $T$-vector of $N\times N$ matrices called \verb"Gvector", then the discretised inflated generator $\mathbf{G}^D$ (denoted \verb"Gdynamic" below) is obtained with: 
%\begin{minted}{julia}
%    using Statistics
%    Gdynamic = mean(Gvector)
%\end{minted}

\subsection{Discretisation of the inflated generator}
\label{sec:discinfgen}
We begin by noting that we can rewrite \eqref{infgenasgen} as
\begin{equation}
    \label{infgen}
\mathcal{G}F(t,x)=G_tF(t,x)+\frac{a^2}{2}\partial_t^2F(t,x).
\end{equation}
Each $G_t$ is estimated as in the previous subsection as $\mathbf{G}^{t_l}$, $l=0,\ldots,T$.
For the second term on the RHS of \eqref{infgen}, we estimate the one-dimensional Laplacian on $[0,\tau]$ by second-order central differencing.
The separation of time and space in the matrix constructions below follows ideas in \cite{AFK24}.
Let $h=t_{l+1}-t_{l}=\tau/T$ be the discrete spacing on the one-dimensional temporal grid in $[0,\tau]$ and $a$ be the temporal diffusion strength parameter.
The indexing we use for the discretised manifold $\mathbb{M}$ is:  the indices of all cubes in $\{t_0\}\times M$, the indices of all cubes in $\{t_1\}\times M$, and so on, ending with the indices of all cubes in $\{t_T\}\times M$.
Including the discretised first components of \eqref{infgen} (namely the $\mathbf{G}^t$) across time and suitably including the discretised second component of \eqref{infgen} (namely a central difference matrix for one-dimensional diffusion) across time, we obtain discretised representations $\mathbf{G}_{\mathrm{spacetime}}$ and $(a^2/2)\mathbf{L}_{\mathrm{spacetime}}$ as below.
Each of these matrices is $TN\times TN$ in size, with $T^2$ blocks that are each of size $N\times N$:
$$\mathbf{G}_{\mathrm{spacetime}}:=\begin{pmatrix}
    \mathbf{G}^{t_0} & 0 & 0&  0 &\cdots & 0  \\
    0&\mathbf{G}^{t_1}& 0 &  0 &\cdots & 0\\
    0 & 0 & \mathbf{G}^{t_2} & 0  & & 0 \\
    0 &0  & 0& \ddots&  &\vdots & \\
        \vdots& & & \ddots&  &\vdots & \\
    0& 0 & \cdots & 0 & 0 &\mathbf{G}^{t_T}
\end{pmatrix}$$
$$\mathbf{L}_{\mathrm{spacetime}}:=\begin{pmatrix}
    -I_N & I_N & 0&  0 &\cdots & 0  \\
    I_N&-2I_N&I_N &  0 &\cdots & 0\\
    0 & I_N & -2I_N & I_N  & & \vdots \\
    0 &0  & I_N & \ddots& \ddots &0  \\
        \vdots&\vdots & & \ddots& -2I_N &I_N  \\
    0& 0 & \cdots & 0 & I_N &-I_N
\end{pmatrix}/h^2.
$$
As per \eqref{infgen} these two matrices are combined to form the discretised inflated generator \begin{equation}
    \label{discinfgen}
\mathbf{G}_a = \mathbf{G}_{\mathrm{spacetime}} + (a^2/2)\mathbf{L}_{\mathrm{spacetime}}.
\end{equation}
If the individual generator matrices $\mathbf{G}^{t_i}$, $i=1,\ldots,T$ are stored in a $T$-vector of $N\times N$ matrices called \verb"Gvector", the matrix $\mathbf{G}_a$ (denoted \verb"Ga" below) is obtained in julia code with:
\begin{minted}[autogobble,fontsize=\small]{julia}
using LinearAlgebra, SparseArrays
Gspacetime = blockdiag(Gvector...)   #block-diagonal array of spatial generators 
Ltemporal = Tridiagonal(ones(T-1), -2*ones(T), ones(T-1))   #1D central difference
Ltemporal[1,1] = -1;  Ltemporal[T,T] = -1   #handle endpoints
Lspacetime = kron(Ltemporal/h^2, I(N))   #kronecker product with N×N identity
Ga = Gspacetime + (a^2/2)*Lspacetime   #combine generator and diffusion components
\end{minted}

%\gf{Mention something about this Kronecker of one-dim temporal Laplacian in \cite{AFK24} by .... }

\subsubsection{Selecting the temporal diffusion strength $a$} \label{ParamSelect}

%Obtaining viable results from the inflated generator method depends on the appropriate assignment of values to the diffusion parameter $a$. We now briefly describe how we choose these values for this study.\
%To compute $\epsilon$, at each time step $t_i$ calculate the median (or mean) of $\| v \left( t_i,x \right) \|$, the velocity magnitudes over all $x \in M$, then average these median (or mean) values over all $T+1$ time steps, and scale this result by a small factor (0.05 was found to be sufficient).
%Selecting a value for $a$ is a little less clear cut by comparison. 
If $a$ is too small, the lack of significant temporal diffusion means that distinct time fibres of a spacetime eigenfunction will be poorly linked together, leading to potentially rapid variation between fibres, which is at odds with our search for quasi-stationary features.
Moreover, temporal eigenvalues will populate the dominant part of the inflated generator spectrum, which is unhelpful for distinguishing spatial features.
If $a$ is too large, then only very slight temporal variation between time fibres in a spacetime eigenfunction will be possible, as small variations will incur high diffusive mixing in the time coordinate due to strong temporal diffusion.
This forced almost-constancy in time may mask quasi-stationary almost-invariant behaviour present in the dynamics. 
In fact, as $a\to\infty$, one expects the eigenspectrum of $\mathcal{G}$ corresponding to non-temporal eignvalues to approach the eigenspectrum of the averaged generator $\bar{G}$,  introduced at the end of Section \ref{sec:2}.
%$\frac{1}{\tau}\int_0^\tau G_t\ dt$, which, due to linearity of the $G_t$ is simply the generator of the time-averaged vector field $\frac{1}{\tau}\int_0^\tau v(\cdot,t)\ dt$. 
\begin{comment}
Selecting a value for $a$ is a little less clear cut by comparison. The value of this parameter should not be too small, or else the temporal diffusion encoded within the inflated generator reduces to the point where a spatial eigenvector at a particular time step is distinct from the eigenvector at all other time steps; and producing a consistent material object over the full time interval from these eigenvectors will prove almost impossible. Whereas if $a$ is too large, the spatial eigenvalues dominate the spectrum to the point where the eigenvectors mirror those obtained from the dynamic generator, and are constant over each time step \cite{FK23}. 
To avoid either of these situations we aim to achieve a healthy medium for $a$, with the best way of doing this involving selecting a value for this parameter which allows for the leading temporal and spatial eigenvalues to be approximately equal. 
\gfadd{As an example of reducing verbosity, the previous sentence can read:} ``To avoid either of these situations we choose $a$ to approximately match the leading temporal and spatial eigenvalues.''
\end{comment}

We now discuss a heuristic for choosing $a$, by matching estimates of leading temporal and spatial eigenvalues.
Recall from the discussion toward the end of Section \ref{sec:qsfamilies} that the leading temporal eigenvalue of the inflated generator is $- a^2\pi^2 /2\tau^2$. 
The real parts of the spatial eigenvalues of $\mathbf{G}_a$ are primarily determined by (i) the spatial diffusion term in (\ref{infgenasgen}) containing $\epsilon$ and (ii) the effective numerical diffusion due to the Ulam scheme.
The former is chosen so that $\epsilon^2/2=0.05\bar{v}\ell$.
For the latter, we estimate effective numerical diffusion, denoted $\epsilon_{\rm eff}$ by $\epsilon_{\rm eff}^2/2=\bar{v}\ell/2$ (see the discussion in Remark 4.6 \cite{FJK13}).
We define a total diffusion $\epsilon_{\rm tot}^2/2:=\epsilon^2/2+\epsilon_{\rm eff}^2/2=0.55\bar{v}\ell$. 
%clearly for small Ulam cubes, the effective numerical diffusion dominates. 
%leading to $\epsilon_{\rm tot}=\sqrt{\bar{v}\ell(\ell/10+1)}$.
Assuming $M$ is a $d$-dimensional hyperrectangle, the leading nontrivial spatial eigenvalue of $(\epsilon_{\rm tot}^2/2)\Delta$ on $M$ is therefore $-\epsilon_{\rm tot}^2\pi^2/2L_{\rm max}^2=-0.55\bar{v}\ell\pi^2/L_{\max}^2$, where $L_{\rm max}$ is the longest side length of $M$. 
As discussed in Section \ref{sec:qsfamilies}, it is desirable to equate diffusion in space and time, so equating these two eigenvalues leads to \begin{equation}
    \label{aeqn}
a=\tau\sqrt{1.1\bar{v}\ell}/L_{\rm max}.
\end{equation}
The value of $a$ can be initialised with this heuristic and adjusted based on subsequent eigenspectrum computations to better match the leading temporal and spatial eigenvalues.
Related discussions on the parameter $a$ for the inflated dynamic Laplacian may be found in Section 4.1 \cite{FK23}, where the estimate of the spatial eigenvalue in \cite{FK23} arises from a Laplace operator, without any considerations of a spatial diffusion parameter $\epsilon$.

\subsubsection{Extracting individual quasi-stationary families with SEBA}

Spectral analysis of dynamical systems to determine almost-invariant sets with transfer operators or generators, and finite-time coherent sets with transfer operators or dynamic Laplacians, sometimes involves several eigenvectors which encode several of these sets in a nontrivial way.
Individual almost-invariant or coherent sets need to be disentangled from the leading eigenvectors and a straightforward way of doing this is with the SEBA algorithm \cite{FRS19}.
The input to SEBA is a collection of eigenvectors and the output is a collection of sparse vectors (whose supports correspond to individual sets)
with approximately the same span as the input eigenvectors. 
SEBA is also very effective in situations where the eigenvectors capture spacetime information, for example separating individual coherent sets \cite{FKS20} and semi-material coherent sets  \cite{FK23,AFK24} across both space and time simultaneously.
In our present setting, if the quasi-stationary family of almost-invariant sets we seek are not easily identifiable from a \emph{single} spacetime eigenvector of the inflated generator, we will apply SEBA to several leading spatial (i.e.\ non-temporal) spacetime eigenvectors to extract several individual quasi-stationary families.

Our overall numerical approach to find quasi-stationary families of almost-invariant sets from a time-dependent velocity field is summarised below.

\begin{algorithm}[Identify quasi-stationary families of almost-invariant sets from vector fields]
\label{alg1}
\, 

\noindent \textit{Inputs:} state space $M$, time-varying velocity field $v$ on $M$, time duration $[0,\tau]$

\noindent \textit{Outputs:} estimates of families of quasi-stationary almost-invariant sets at discrete time points in $[0,\tau]$ 
\begin{enumerate}
   \item Discretise the space space $M$ into cubes $B_1,\ldots,B_N$, and the time interval $[0,\tau]$ into time nodes $t_{l}$, $l=0,\ldots,T$ as described at the beginning of Section \ref{sec:disc}.
   \item Select a value for $\epsilon$ using equation \eqref{epseqn} and compute the discrete generator matrices $\mathbf{G}^{t_l}$ for $l=0,\ldots,T$ as described in Section \ref{sec:gendisc} using \eqref{Geqn}.
   \item Combine the $\mathbf{G}^{t_l}$ to form $\mathbf{G}_a$ as described in Section \ref{sec:discinfgen} using \eqref{discinfgen}, where the temporal strength parameter $a$ is initially selected using \eqref{aeqn}.
   \item Compute the real eigenvalues $0=\Lambda_1>\Lambda_2\ge \cdots$ of the sparse matrix $\mathbf{G}_a^*$ that are closest to 0, along with their corresponding spacetime eigenvectors $\mathbf{F}_1^*,\mathbf{F}_2^*,\ldots$.  Note that $\mathbf{F}_1^*$ should be a constant vector $\mathbf{1}$.
   \item Remove the constant-in-space temporal eigenvectors $\mathbf{F}_k^*$ by inspection or numerical interrogation; see the discussion at the end of Section \ref{sec:qsfamilies} for their functional form. The discretised temporal eigenvectors will be constant on blocks of indices $lN+1,\ldots,(l+1)N$, $l=0,\ldots,T$.
   \item If the desired quasi-stationary family of almost-invariant sets is not clearly apparent in any single eigenvector, but there are signatures of the family (i.e.\ large positive and negative values indicate relevant spatial structures) in an eigenvector $\mathbf{F}_{K}^*$ for some $K\ge 2$, apply SEBA to the $K'$ real-valued \textit{spatial} eigenvectors with index $K$ and below to obtain spacetime SEBA vectors $\mathbf{S}_1,\ldots,\mathbf{S}_{K'}$, with $K'\le K$.  
   One of the $\mathbf{S}_k$, $k=1,\ldots,K'$ should be supported on the desired quasi-stationary family.
\end{enumerate}
\end{algorithm}
\section{An idealised flow with sharp variation in quasi-stationary structures} \label{IdealisedDGFlow}

% All relevant results should be included, documenting the evolution of our research to a final result. At the moment, we have two systems of interest:

In this section we illustrate our constructions on an idealised fluid flow.
We use the ``switching double gyre'' flow considered in Section 6.1 \cite{AFK24}, which is a modification of the standard double gyre \cite{Shaddenetal05}.
\begin{figure}[htb]
      \centering
      \includegraphics[width=\textwidth]{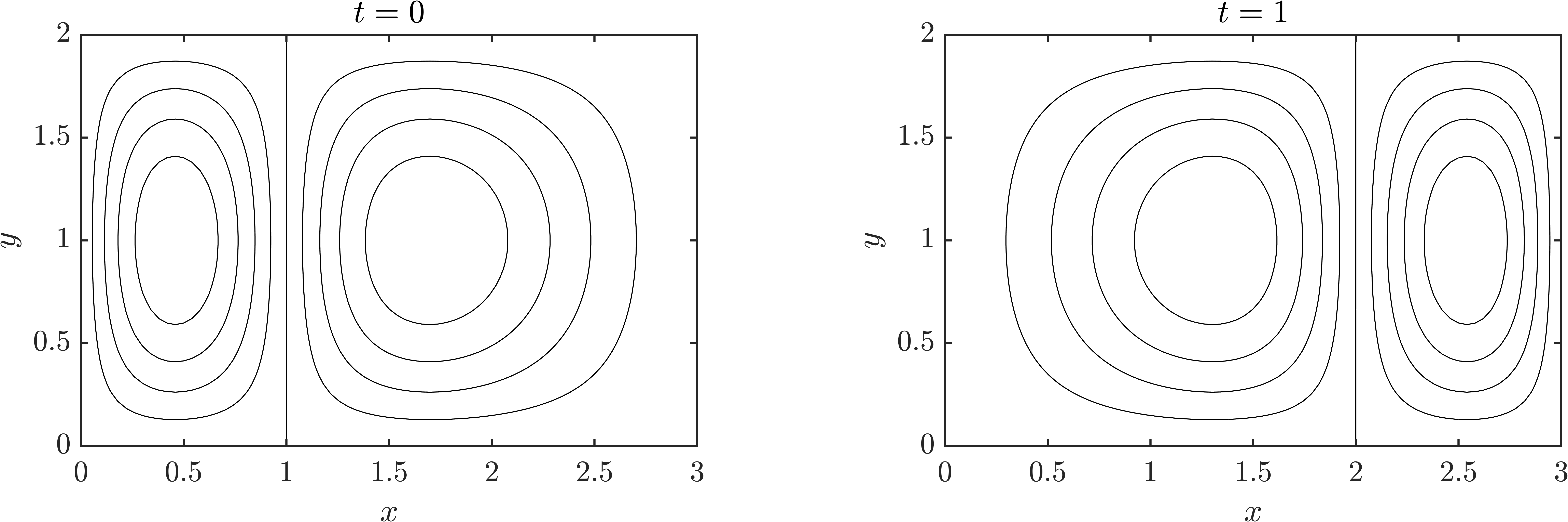}
      \caption{Some level sets of the stream function for the switching double gyre flow taken at $t = 0$ (left) and $t = 1$ (right).}
      \label{fig:DG_Streamplot}
\end{figure}
The domain is $M = \left[0,3\right] \times \left[0,2\right]$,
%and the flow time runs from 0 to $\tau = 1$.
%For the definition of the flow, we refer to \cite{AFK24}.
and the main almost-invariant features of the flow are two gyres of unequal size. 
Initially the left gyre $A_0^{\rm left}$ is half the size of the right gyre $A_0^{\rm right}$, and as the flow progresses from $t=0$ to $t=1$, the sizes change sharply at times between approximately $t=0.4$ and $t=0.6$, so that at the end of the flow duration the relative sizes have switched over:  $A_1^{\rm left}$ is double the size of $A_1^{\rm right}$;  see Figure \ref{fig:DG_Streamplot} for streamline plots at the initial and final times.

The paper \cite{AFK24} studied the coherence properties of the gyres, attempting to track the evolution of the gyres in space as their volume changed dramatically.
In the present paper we are instead concerned with determining the existence of a quasi-stationary family of almost-invariant objects through time.
In contrast to coherent sets, which may freely move about in space while remaining coherent, a quasi-stationary family of almost-invariant sets are approximately fixed in space. 
Thus, the dynamical question we pose is different, and particularly relevant for e.g.\ atmospheric blocks, which remain approximately fixed in space.

\begin{figure}[htb]
      \centering
      \includegraphics[width=0.5\textwidth]{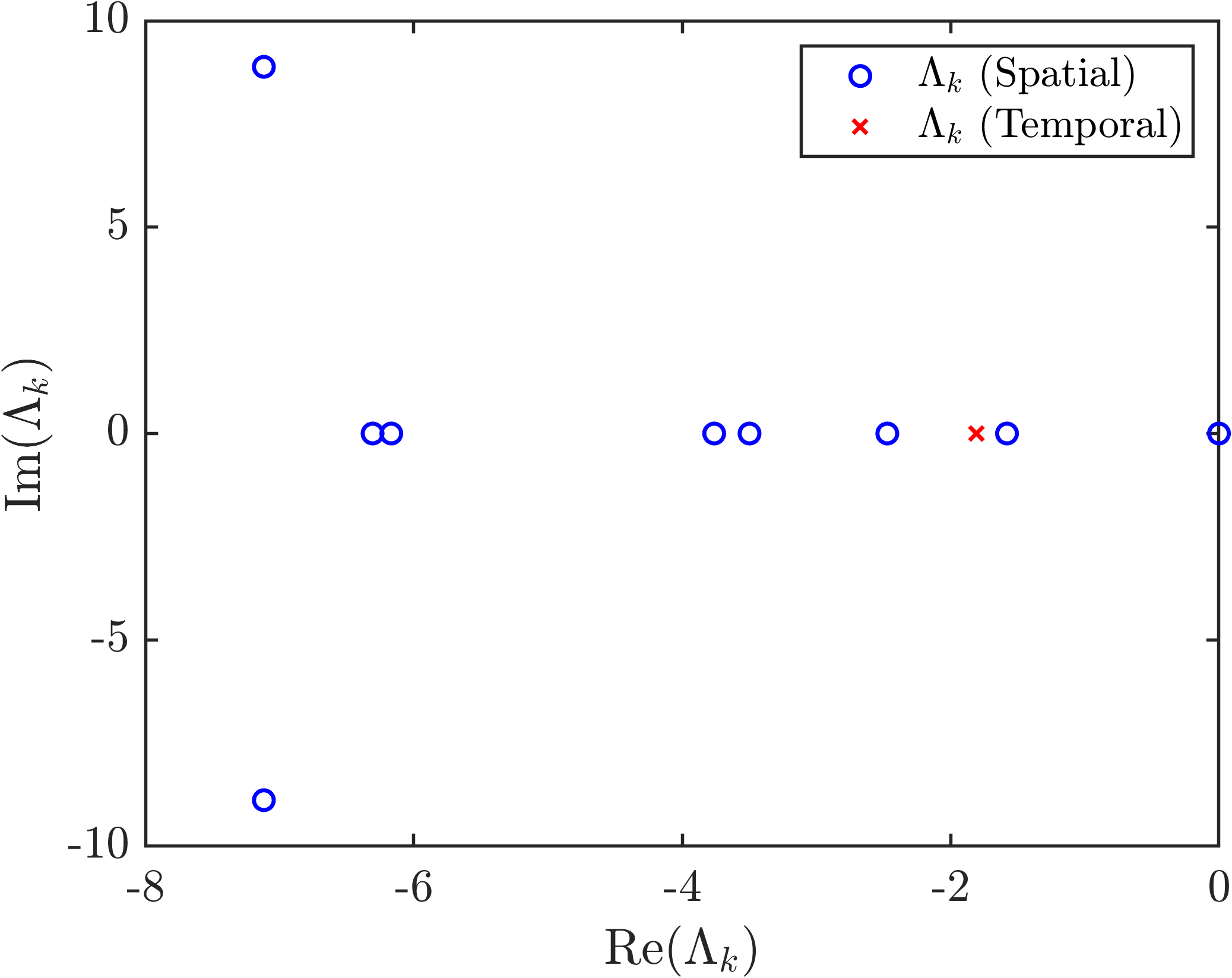}
      \caption{The leading ten eigenvalues of the inflated generator approximation $\mathbf{G}_a$ for the switching double gyre flow.}
      \label{fig:DG_Eigvals}
\end{figure}

\begin{figure}[htb]
    \centering
    \includegraphics[width=\textwidth]{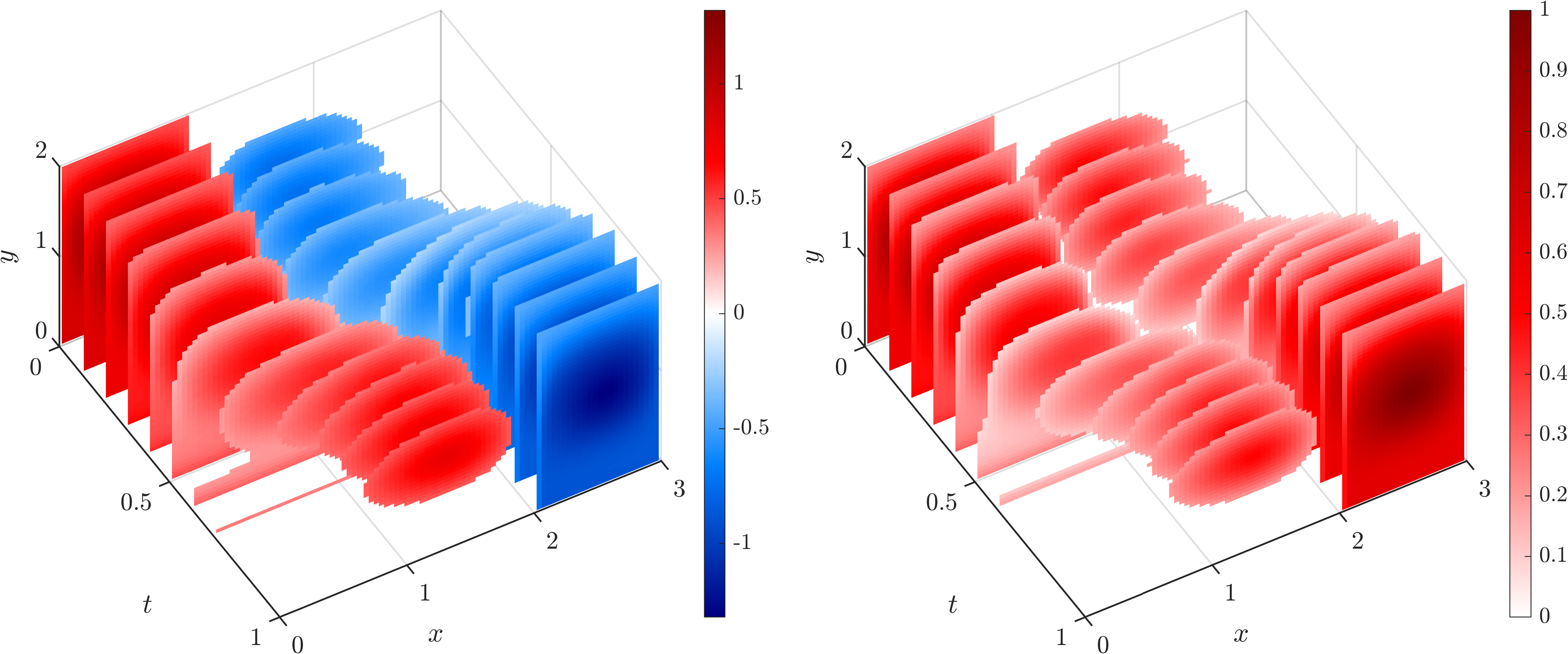}
    \caption{Inflated generator construction of quasi-invariant families of almost-invariant sets in the switching double gyre. \emph{Left:} Selected time fibres of the second eigenvector $\mathbf{F}^*_2$ of $\mathbf{G}_a^*$.  A cutoff of 0.33 for the magnitude of the eigenvector has been applied. \emph{Right:} Superposition of two SEBA vectors $\mathbf{S}_1$ and $\mathbf{S}_2$. A cutoff of 0.1 for the magnitude of the vector superposition has been applied.}
    \label{fig:DGswitch3d}
\end{figure}

For our generator calculations in this experiment, we follow Algorithm \ref{alg1}, discretising $M$ into a $75 \times 50$ grid of square boxes $B_1,\ldots,B_{3750}$ of side length $\ell=0.04$. 
We discretise the time interval $[0,1]$ into 21 nodes spaced $h=0.05$ time units apart.
Following steps 2 and 3 of Algorithm \ref{alg1}, we calculate $\epsilon \approx 0.2397$ from $\bar{v} \approx 14.3698$ and $\ell$.
We set $a = 0.45$ to better match the leading spatial and temporal eigenvalues, see Figure \ref{fig:DG_Eigvals}.
An initial $a$ value of approximately $0.2651$ was calculated using the heuristic in \eqref{aeqn}.
The sparse $3750\times 3750$ matrices $\mathbf{G}^{t_0},\ldots,\mathbf{G}^{t_{20}}$ are combined to form the sparse $78750\times 78750$ matrix $\mathbf{G}_a$.
Figure \ref{fig:DG_Eigvals} shows the leading ten eigenvalues $0=\Lambda_1>\Lambda_2>\cdots>\Lambda_{10}$;  $\Lambda_3$ is a temporal eigenvalue and is shown as a red cross.
Figure \ref{fig:DGswitch3d}(left) plots the estimates of the second eigenfunction $\mathbf{F}_2^*$ at time instances $0, 0.1, 0.2,\ldots, 1$.
% Old Figure: DoubleGyre_EigvecWithSEBA_3D_LowerCutOff_22Apr24_2.png
% Old Cutoffs: 0.15, 0.01
It will be useful for the reader to refer to Figures \ref{fig:schematic} and \ref{fig:DGswitch3d} during the following discussion.
At any time $t\in [0,1]$, there is no trajectory escape from either rectangular gyre $A_t^{\rm left}$ or $A_t^{\rm right}$ in $M$ under advective from the steady vector field $v(t,\cdot)$. 
The only possibility of trajectory escape from $A_t^{\rm left}$ or $A_t^{\rm right}$ (restricted to $\{t\}\times M$) is by small $\epsilon$-diffusion.
On the other hand, the families $\{A_t^{\rm left}\}_{t\in[0,1]}$ and $\{A_t^{\rm right}\}_{t\in[0,1]}$ undergo a sharp change in spatial position and area for $t\in[0.4,0.6]$.
Therefore if we simply defined $\mathbb{A}^{\rm left}=\cup_{t\in [0,1]} \{t\}\times A_t^{\rm left} \subset\mathbb{M}$ (and defined $\mathbb{A}^{\rm right}$ in a similar fashion), the resulting spacetime set would have a large loss of trajectories due to the $a$-diffusion in the time coordinate, along which the set varies rapidly between times 0.4 and 0.6. 
This na\"{i}ve choice of $\mathbb{A}$ would be optimal in terms of almost-invariance in space, but rather suboptimal in terms of almost-invariance in time.

Turning to our inflated generator results, in Figure \ref{fig:DGswitch3d}(left) we plot time slices of $\mathbf{F}_2^*$ across time.
One can interpret the red sets as time fibres of one candidate almost-invariant set $\mathbb{A}^{(1)}$ and the blue as time fibres of another candidate almost-invariant set $\mathbb{A}^{(2)}$.
Note that parts of $\mathbb{A}^{(1)}$ are cut away (compared to $\mathbb{A}^{\rm left}$) at times from about $t=0.4$ onwards, when the velocity field $v(t,\cdot)$ first undergoes rapid change, moving the separatrix in the vector field to increase the gyre area.
This cutting away is designed to improve the \textit{spacetime} almost-invariance of the spacetime set $\mathbb{A}^{(1)}$, and therefore the \textit{quasi-stationarity} of the family $\{A^{(1)}_t\}_{t\in [0,1]}$.
This comes at the cost of a slight decrease in the spatial almost-invariance of $A_t^{(1)}$ (compared to $A_t^{\rm left}$) for some $t$, but overall our inflated generator construction identifies an improved balance between almost-invariance in space and time.
The same interpretation can be given to $\mathbb{A}^{(2)}$, where some parts are cut away, relative to $\mathbb{A}^{\rm right}$, from time $t=0$ up to about $t=0.6$.

Inserting $\mathbf{F}^*_1$ and $\mathbf{F}^*_2$ into SEBA yields two SEBA vectors  $\mathbf{S}_1, \mathbf{S}_2$, highlighting the individual spacetime sets $\mathbb{A}^{(1)}, \mathbb{A}^{(2)}$, respectively; see Figure \ref{fig:DGswitch3d}(right).
There is a variation in the depth of colour in the SEBA vectors.
The deepest red indicates those times where one is most certain that points at those spacetime locations belong to the quasi-stationary family of almost-invariant sets.
We see in $\mathbb{A}^{(1)}$ that the time fibres for $t$ near the initial time $t=0$ consist of points that strongly belong to this family.
This is because these points belong to a small initial gyre, which will very likely remain in the family of almost-invariant sets as this family approximately expands about it in space.
On the other hand, points in the final time fibres of $\mathbb{A}^{(1)}$ (those near $t=1$) are paler and less certain to remain in the family;  points nearer to the gyre ``core" are deepest red and most certain to remain.

% we provide three-dimensional visualisations of $V_2$ and the superposition of two SEBA vectors $\Sigma_1$ and $\Sigma_2$ produced using the first two eigenvectors of the inflated generator. The inflated generator method identifies quasi--stationary flow objects as the evolutions of said objects as we progress through time over the flow window of interest. The results more accurately illustrate how key patterns of the flow system evolve, and in this case they do so in the precise fashion described at the beginning of this section.

\section{Blocking over Europe}
\label{sec:block}
One important practical application of this work involves identifying the shape and tracing the movement of atmospheric blocking events. We will identify blocking events by applying the inflated generator method to atmospheric wind velocity data, and then compare these results to those obtained from three well-documented block detection methods. We also compare these results to wind streamlines, \new{in particular to the anticyclonic wind circulation patterns expected within blocks \cite{AnticycloneRotation}, and to observed surface temperature anomaly fields to illustrate how blocking events may be correlated with extreme weather.}

We use ERA5 atmospheric wind velocity, geopotential height, potential vorticity \new{and temperature} reanalysis data obtained from the Copernicus Climate Change Service Climate Data Store \cite{ERA5,ERA5_singlelevels}; over the spatial region $M_{\mathrm{full}} = \left[ 20^{\circ} \mathrm{W} , 60^{\circ} \mathrm{E} \right] \times \left[ 15^{\circ} \mathrm{N} , 90^{\circ} \mathrm{N} \right]$ and the time period ranging from 00:00 01/01/1979 to 18:00 31/12/2019 (UTC) inclusive. \new{We require potential vorticity, geopotential height and temperature data over this multidecadal time span to compute anomalies of these quantities.
In our generator computations, we only use wind velocity data from 00:00 01/07/2003 to 18:00 31/08/2003 (UTC).} The potential vorticity data is taken on all available vertical pressure levels ranging from 150 to 500 hPa (inclusive), while the wind velocity and geopotential height data are taken at the 500 hPa pressure level only. \new{Temperature data is taken at the surface of the Earth. For consistency with the vast majority of the existing blocking literature, we represent the blocks in two dimensions at the 500 hPa vertical pressure level.} The spatial resolution of this data is $0.25^{\circ}$ in both longitude and latitude, and the temporal resolution for the data is 6 hours. 
%\gf{Aleks to briefly mention interpolation here and the use of Interpolations.jl (and cite it).} \aleks{Here is the text:} 
We linearly interpolate velocities between the spatial data points %we  \verb"linear_interpolation" convenience constructor featured 
with julia's \verb"Interpolations.jl" package \cite{Interpolations}. %A convenience constructor for higher order cubic spline interpolation is also available within this package, however after some experimentation linear interpolation was considered sufficient for our calculations. 
%\gf{how did you decided on attributed author names?  You included version number for this but not for quadgk?} \aleks{The GitHub repositories for each of these packages contain Citation files (in either .bib or .eff format) which show you how to cite the packages and list the authors to cite. Here are the links: https://github.com/JuliaMath/Interpolations.jl/blob/master/CITATION.cff ; https://github.com/JuliaMath/QuadGK.jl/blob/master/CITATION.bib. The QuadGK.jl citation does not require you to list the version number (refer to the link), but I can include it: it's v2.9.4.} 
We use the full data set for calculations involving geopotential height, potential vorticity,  \new{and temperature;} but for our inflated generator calculations we do not use the entire spatial domain nor the full spatial resolution available.

\subsection{Existing Block Detection Methods} \label{BlockIdentification}

We begin by briefly describing three existing geophysical block detection methods. The reader is referred to the cited literature for more information on each method. \new{We also describe how surface temperature anomalies are computed in this section, though this is not a block detection proxy.}

\begin{comment}
We start by describing in detail three of the most commonly used methods from the literature to identify blocking events. Blocks obtained from these methods will then be compared against the blocks obtained from the generator methods.
\end{comment}

\paragraph{Geopotential height gradients.} Gradients of $\mathrm{Z}_{500}$, %\gf{Is $\zeta$ standard notation?} \aleks{No, $Z_{500}$ is. I'll adjust below accordingly.}
which denotes geopotential height taken at a 500 hPa pressure level. While this method was originally used to identify blocked longitude levels through changes in the $\mathrm{Z}_{500}$ gradient over a single line of latitude \cite{Lejenas83,Tibaldi90}, we instead consider a ``two--dimensional" version of this method \cite{Scherrer06,Buehler11,Pinheiro19,Thomas21} used to identify blocked regions at a particular moment in time. Let $\lambda$ and $\phi$ represent  longitude and latitude respectively. Using ``northern" and ``southern" $\mathrm{Z}_{500}$ gradients, as defined in \cite{Scherrer06} and denoted here by $\gamma_N$ and $\gamma_S$ respectively, %\aleks{The usual notation for these gradients simply involves defining these as ``GHGN" and ``GHGS" (respectively).}
a closed set $\Gamma_t \subset M$ at an instantaneous moment in time $t$, where
\begin{equation*}
    \Gamma_t = \left\{ \left( \lambda,\phi \right) \in M \ : \ \gamma_N \left( t,\lambda,\phi \right) < -10, \ \gamma_S \left( t,\lambda,\phi \right) > 0 \right\}
\end{equation*}
is classified as a blocked region. 

\paragraph{Geopotential height anomalies.} This method also utilises $\mathrm{Z}_{500}$ through the computation of the long--term \textit{anomalies} of this quantity, which we denote by $\mathrm{Z}_{500}^{*}$. 
%\aleks{A separate symbol for geopotential height/potential vorticity anomalies is seldom used in the literature; the authors often write out ``$Z_{500}$ anomalies" each time.}. 
We compute $\mathrm{Z}_{500}^{*}$ using the procedure detailed in \cite{Sausen95}, and identify blocking events as closed contours of $\mathrm{Z}_{500}^{*}$ that are above a threshold of +100 geopotential metres (gpm) \cite{Dole83,Carrera04,Pinheiro19,Thomas21}.

\paragraph{Vertically averaged potential vorticity anomalies.}  Anomalies of potential vorticity averaged over a range of vertical pressure levels, $\overline{\mathrm{VPV}}$, 
%\gf{why do we use $\omega$ for PV? $\omega$ is typically for vorticity.  Why do we need the subscript VA if you also use a tilde?} \aleks{``VA" stands for ``vertically averaged". A Greek letter isn't used for the potential vorticity usually, rather it is denoted by VPV (Vertically (averaged) Potential Vorticity) or VAPV.}
which we denote by $\overline{\mathrm{VPV}}^{*}$.
%\aleks{$\mathrm{VPV}^*$ is often used to denote these anomalies.} 
%\gf{Maybe it is more meaningful to use $\overline{\mathrm{VPV}^*}$ and $\mathrm{Z}_{500}^*$ (for gpht anomalies), so the star is the anomaly and the bar indicates an average?} 
These anomalies are computed using the methodology described in \cite{Pinheiro19}, with the omission of the Fourier transform harmonics step used to smoothen the ``Long Term Daily Means" of $\overline{\mathrm{VPV}}$. We computed $\overline{\mathrm{VPV}}^{*}$ without this extra step and found the results to be very similar to those in \cite{Thomas21}.
%, we were satisfied with the results produced and we deemed this step to be unnecessary. 
We then identify blocks similarly to the method of the previous paragraph, except we seek contours of $\overline{\mathrm{VPV}}^{*}$ below a threshold of -1.2 potential vorticity units \cite{Schwierz04,Pinheiro19,Thomas21}.
\vspace*{.1cm}

Table \ref{tab:BlockingEvents} lists two blocking events occurring over Europe between July and August 2003 identified in \cite{Thomas21} using at least one of these three methods and a ``Ground Truth Dataset" used to verify all identified blocks. Rows 2, 3 and 4 of Figures \ref{fig:EuroBlock_East_MultiFig} and \ref{fig:EuroBlock_West_MultiFig} show instantaneous snapshots of regions $\Gamma_t$ satisfying the $\mathrm{Z}_{500}$ gradient criteria and scalar fields of $\mathrm{Z}_{500}^{*}$ and $\overline{\mathrm{VPV}}^{*}$ respectively, which visually identify blocks. The blocking criteria described above should  be maintained for a fixed amount of time (usually between four \cite{Tibaldi90,Pelly03,Ehstand21} to ten \cite{Rex50II,Dole83} days) for a region to be considered blocked, but in practice this condition is often relaxed. 
%but we will not adhere to this extra criterion too strictly as one can intuitively develop an image of a consistently blocked region through qualitative observation of these individual instances of each set/scalar field.

\begin{table}[h!]
    \centering
    \begin{tabular}{|c|c|c|c|c|}
    \hline
       \textbf{Block} & \textbf{Date Range} & \textbf{Long. Range} & \textbf{Lat. Range} & \textbf{Identified Using} \\
    \hline
       East & 29/07/03 -- 03/08/03 & $25^{\circ}$ E -- $55^{\circ}$ E & $55^{\circ}$ N -- $70^{\circ}$ N & $\Gamma_{t}$, $\mathrm{Z}_{500}^{*}$, $\overline{\mathrm{VPV}}^{*}$ \\
    \hline
    West & 02/08/03 -- 12/08/03 & $10^{\circ}$ W -- $15^{\circ}$ E & $40^{\circ}$ N -- $60^{\circ}$ N & $\mathrm{Z}_{500}^{*}$ \\
    \hline
    \end{tabular}
    \caption{The ``east" and ``west" blocking events as identified in \cite{Thomas21}. The longitudinal and latitudinal ranges of these blocks are estimated from Figure 5 of \cite{Thomas21}.}
    \label{tab:BlockingEvents}
\end{table}

%\new{\paragraph{Surface temperature anomalies.}  
\new{Atmospheric blocks can sometimes be associated with heatwaves, as is the case for the blocks listed in Table \ref{tab:BlockingEvents} \cite{Ehstand21,Thomas21}. 
To check if there is a correlation between the East and West blocks and anomalous surface temperatures, we calculate a standard \cite{AnomalyDefn,AnomalyComp} surface temperature anomaly field $T_{\mathrm{surf}}^{*}$ as follows. Long-term daily temperature means are computed using ERA5 surface temperatures ($T_{\mathrm{surf}}$) \cite{ERA5_singlelevels}, recorded at 00:00 (UTC) between 01/01/1979 and 31/12/2019 to maintain consistency with the anomaly calculations for geopotential height and potential vorticity. Surface temperature anomalies are then obtained by subtracting these long-term daily means from the surface temperatures recorded at 00:00 (UTC) between 01/07/2003 and 31/08/2003.}

\subsection{Results of block identification using the inflated generator} \label{EuroBlock}

We wish to apply the inflated generator method to atmospheric velocity data to identify the onset, evolution, and dissipation of the two blocking events listed in Table \ref{tab:BlockingEvents}, documented \cite{Ehstand21,Thomas21} as having played a key role in the 2003 European heatwaves.
We will take $M$ to be a fixed rectangle restricted to a surface corresponding to a pressure level of 500 hPa \cite{ERA5} and $v$ to be the time-varying wind field on this surface.
The vector field $v$ need not have zero normal velocity at the boundary of $M$.
Following the approach in section \ref{sec:gendisc} we discretise the 500 hPa pressure level into coordinates in degrees longitude and latitude, which results in  rectangular boxes $B_1,\ldots,B_N$ with areas that depend on the latitude.
We denote the longitudinal and latitudinal extent (in metres) of box $B_i$ by $\ell_{i,{\rm lon}}$ and $\ell_{i,{\rm lat}}$, respectively.
Equation \eqref{Geqn} becomes:
\begin{equation}
    \label{Geqngeneral2d}
\mathbf{G}^{t_l}_{ij}=\begin{cases}
    \frac{1}{m(B_j)}\left(\int_{B_i\cap B_j}(\max\{v(t_l,x)\cdot n_{ij},0\})\ dm_{1}(x)\right) + \epsilon^2/(2\ell_{i,{\rm lon}}^2),&\quad i\neq j, \mbox{$n_{ij}$ is longitudinal}\\
    \frac{1}{m(B_j)}\left(\int_{B_i\cap B_j}(\max\{v(t_l,x)\cdot n_{ij},0\})\ dm_{1}(x)\right) + \epsilon^2/(2\ell_{i,{\rm lat}}^2),&\quad i\neq j, \mbox{$n_{ij}$ is latitudinal}\\
    \sum_{j\neq i}-\frac{m(B_j)}{m(B_i)}\mathbf{G}^{t_l}_{ij},&\quad\mbox{otherwise.}
\end{cases}
\end{equation}
\new{For the purposes of choosing an initial value of $a$,} because our box sizes vary with latitude, we denote the median box side length (computed over longitudinal and latitudinal directions) by $\bar{\ell}$.
%\new{This produces a uniform diffusion  the same noise level is applied uniformly in space.}
%\new{GF TO UPDATE THIS. This variation in box size does not introduce any bias in our results, as the area of each box is accounted for in our flux computations and the same noise level is applied uniformly in space.}
Analogously to the choice in Section \ref{sec:gendisc} the small spatial diffusion parameter \new{$\epsilon_{\rm tot}$ is set to $\epsilon_{\rm tot}=\sqrt{1.1\bar{v}\bar{\ell}}$.}
% \new{GF TO ADD EPS-TOT.}
% The former is chosen so that $\epsilon^2/2=0.05\bar{v}\ell$.
% For the latter, we estimate effective numerical diffusion, denoted $\epsilon_{\rm eff}$ by $\epsilon_{\rm eff}^2/2=\bar{v}\ell/2$ (see the discussion in Remark 4.6 \cite{FJK13}).
% We define a total diffusion $\epsilon_{\rm tot}^2/2:=\epsilon^2/2+\epsilon_{\rm eff}^2/2=0.55\bar{v}\ell$. 
\new{The initial guess for the temporal diffusion strength parameter $a$ is correspondingly modified from the discussion in Section \ref{ParamSelect}; 
 we set $a=\tau\sqrt{1.1\bar{v}\bar{\ell}}/L_{\rm max}$.}
%, one occurring further West and centred around $5^{\circ} \mathrm{E}$ (referred to as the ``Western Block"); and one occurring further East and centred around $40^{\circ} \mathrm{E}$ (referred to as the ``Eastern Block"). 
%The results obtained will then be compared against the other block detection techniques described in the previous subsection to see how well the inflated generator can identify these blocking events. In addition, we will also compare the inflated generator results with the wind speed scalar field and its associated streamlines.

\subsubsection{East Block (July-August 2003)} \label{EuroBlockEast}

The first blocking event we consider occurred in Eastern Europe between 00:00 29/07/03 to 00:00 03/08/03 (UTC). To identify this block, we let $M = \left[ 15^{\circ} \mathrm{E} , 60^{\circ} \mathrm{E} \right] \times \left[ 30^{\circ} \mathrm{N} , 75^{\circ} \mathrm{N} \right]$. The temporal domain includes the full duration of the block with a three day extension applied to either side of this time span; giving an interval ranging from 00:00 26/07/03 to 00:00 06/08/03 (UTC) (inclusive) of length $\tau = 11$ days. We use the full temporal resolution of 6 hours, and a spatial resolution of $1^{\circ}$ in both latitude and longitude; providing a spatial grid of $45 \times 45$ boxes.

\begin{comment}

We start once more with the dynamic generator method, with the results obtained from this method shown in Figure \ref{fig:EuroBlock_East_DynGen}. On display in this Figure is the first non--trivial eigenvector of the generator $v_2$, along with $S_1$, the first of two SEBA vectors produced using the first two eigenvectors of the dynamic generator. As desired, the minimal value of each of the two SEBA vectors $S_1$ and $S_2$ is zero. Unlike in the Western block case, the dynamic generator results have been able to successfully identify the Eastern block, as the elliptic red shape shown in both panels takes the expected shape of the block described in Table \ref{BlockEventList} and in \cite{Thomas21}. However, like in Section \ref{IdealisedDGFlow}, it must once again be stressed that the outcome of this method by design gives the impression that the coherent object identified is not only present but also maintains this shape and size over the course of the entire time interval, which is certainly not the case.

\begin{figure}[h]
      \centering
      \includegraphics[width=0.9\textwidth]{EuroBlock_East_GeneratorEigVec_19Mar24.png}
      \caption{The leading non-trivial eigenvector $v_2$ and the first SEBA vector $S_1$ obtained from the dynamic generator for the European atmospheric flow (Eastern region) over the date range 26/07/03-06/08/03.}
      \label{fig:EuroBlock_East_DynGen}
\end{figure}

\end{comment}

\begin{figure}[hbt!]
      \centering
      \includegraphics[width=0.5\textwidth]{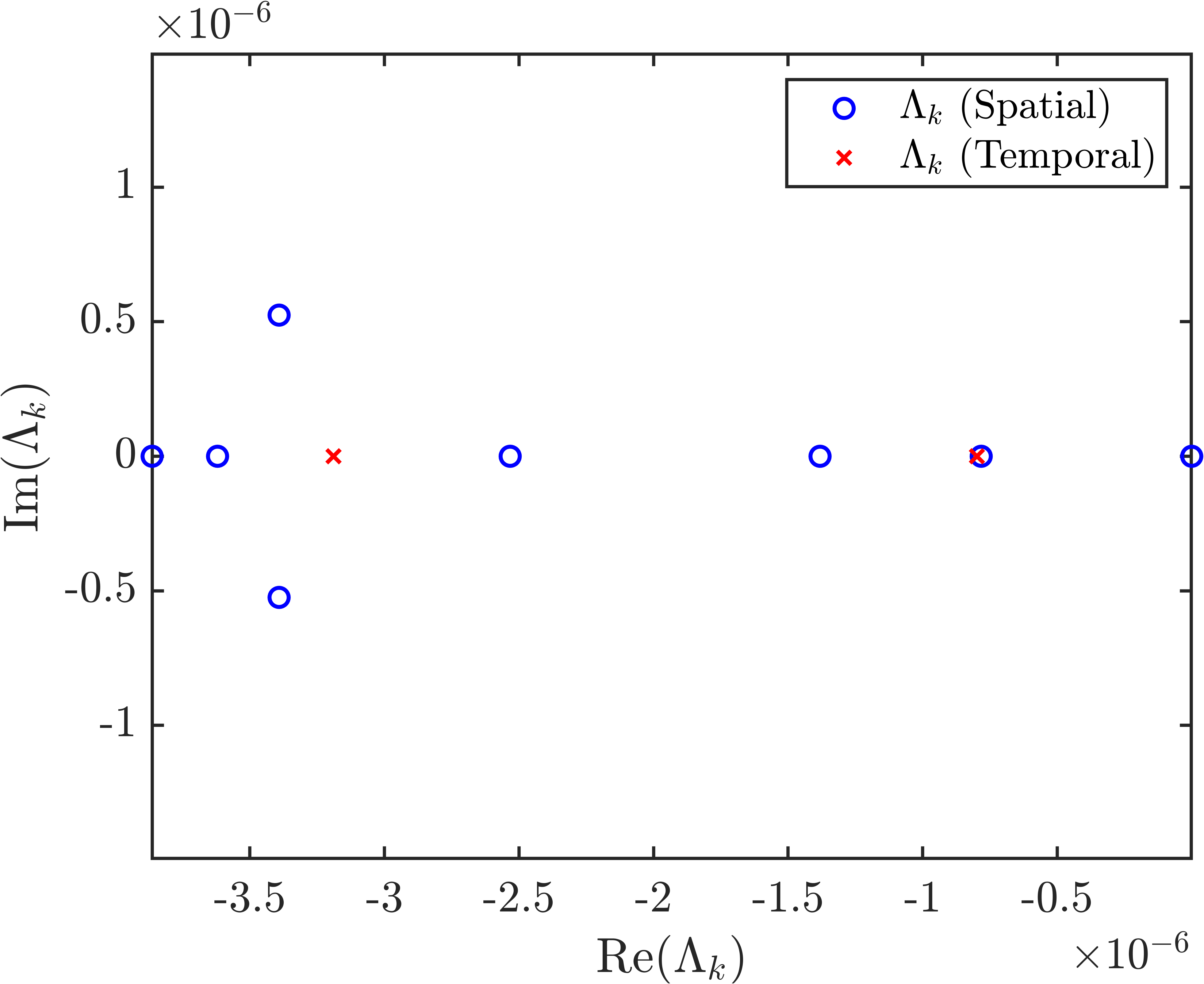}
      \caption{The leading ten eigenvalues $\Lambda_k$ (ordered according to their real parts in decreasing order) of the inflated generator calculated for the European atmospheric flow (eastern region) over the date range 26/07/03-06/08/03, plotted in the complex plane. Spatial and temporal eigenvalues are distinguished using separate plot markers and colours.}
      \label{fig:EuroBlock_East_Eigvals}
\end{figure}

\begin{figure}[hbt!]
      \centering
\includegraphics[width=\textwidth]{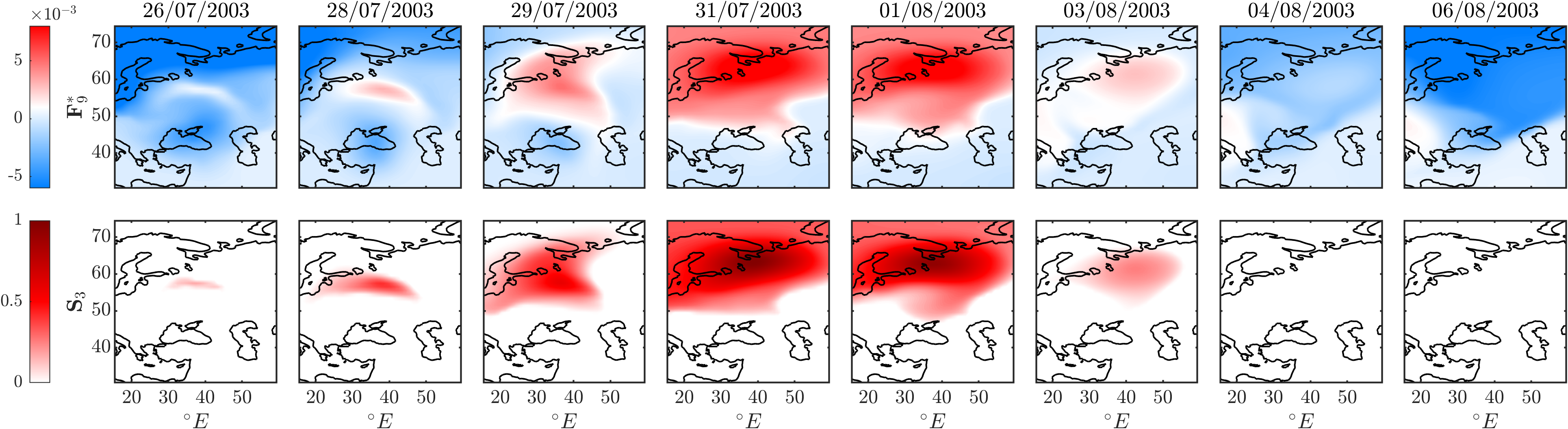}
      \caption{Temporal evolution of the spatial eigenvector $\mathbf{F}_{9}^{*}$ and the SEBA vector $\mathbf{S}_3$ of the inflated generator for the European atmospheric flow (eastern region) at select dates taken over the date range 26/07/03-06/08/03. The time of day is 00:00 (UTC) for all of the above dates.}
      \label{fig:EuroBlock_East_InfGen}
\end{figure}

Following the instructions described in Section \ref{ParamSelect}, we calculate $\epsilon \approx 319.2460 \ \mathrm{ms}^{-1/2}$, given a median box side length $\bar{\ell} \approx 103618 \ \mathrm{m}$ and median speed $\bar{v} \approx 9.8360 \ \mathrm{ms^{-1}}$, and let $a = 0.0032$;  our heuristic yielded $a=0.0023$. The leading ten eigenvalues for the inflated generator are shown in Figure \ref{fig:EuroBlock_East_Eigvals}.
\new{Eigenvectors $\mathbf{F}_3^*$ and $\mathbf{F}_6^*$ are real valued and temporal, while $\mathbf{F}_7^*$ and $\mathbf{F}_8^*$ are complex valued. 
All of the other leading ten eigenvectors are real valued and spatial, and we use these to identify quasi-stationary families of almost-invariant sets.}
We then visually inspect each of the real-valued spatial eigenvectors individually to identify which of these generates a quasi-stationary object that best resembles the east block. The deep red elliptical objects identified from the ninth eigenvector $\mathbf{F}_{9}^{*}$ shown in Figure \ref{fig:EuroBlock_East_InfGen} (upper row) were found to provide the best visualisation of this block. By comparison, in results not shown, $\mathbf{F}_2^{*}$ identifies a combination of the east block and another quasi-stationary object positioned further east; while $\mathbf{F}_4^{*}$ and $\mathbf{F}_5^{*}$ each identify a cyclonic object nested underneath the east block. The latter object is the cut-off low expected to appear underneath a block \cite{Shutts83}. We then apply the SEBA algorithm to all real-valued spatial eigenvectors from $\mathbf{F}_{1}^{*}$ up to and including $\mathbf{F}_{9}^{*}$ (namely $\mathbf{F}_{1}^{*},\mathbf{F}_{2}^{*},\mathbf{F}_{4}^{*},\mathbf{F}_{5}^{*},\mathbf{F}_{9}^{*}$), and then inspect the five SEBA vectors to identify which of these gives the best visualisation of the east block. 
The third SEBA vector $\mathbf{S}_3$, shown in Figure \ref{fig:EuroBlock_East_InfGen} (lower row) isolates this block.
Figure \ref{fig:EuroBlock_East_InfGen} (lower row) shows that the east block gradually begins to form around July 26, with the block taking its expected shape just after July 29. 
The  block is maintained until August 3, by which time the block has reduced to a smaller elliptical object. 
No trace of the block can be identified from $\mathbf{F}_{9}^{*}$ or $\mathbf{S}_3$ beyond August 3.

\begin{comment}
\textcolor{Cerulean}{Figure \ref{fig:EuroBlock_East_3DayExtension_iDL_Results} shows the results used to identify the Eastern block for this time interval. The block itself does not appear to be a dominant quasi--stationary object within this time interval, seeing as we need to use the tenth dynamic generator eigenvector $V_{10}$ to identify it. We use six eigenvectors $V_1, V_3, V_4, V_6, V_7$ and $V_{10}$ in the SEBA algorithm, excluding all temporal and complex spatial eigenvectors, and obtain six SEBA vectors with minimum values $0, 0, -0.0077, -0.0675, -0.1731$ and $-0.3456$ respectively. While some of these minima are a bit small, including all of these eigenvectors is necessary to obtain the block through SEBA, which we are able to do through the fourth SEBA vector $\Sigma_4$. There is no trace of any blocking activity at the start of the domain (as expected), with the block gradually building up over time before taking full formation on July 30, one day earlier than in the non-extended case.} \aleks{Figures will be updated later on to include August 29 results.} \textcolor{Cerulean}{We obtain the elliptical shape expected from this block between July 30 and August 3, by which time the block has reduced to a small elliptical object of low intensity. No trace of the block can be identified from $V_{10}$ or $\Sigma_4$ beyond August 3. There is a noticeable contrast between the blocking object identifed for this time domain and the non-extended time domain considered in the previous section.}
\end{comment}

\begin{figure}[hbt!]
    \centering
    \includegraphics[width=\textwidth]{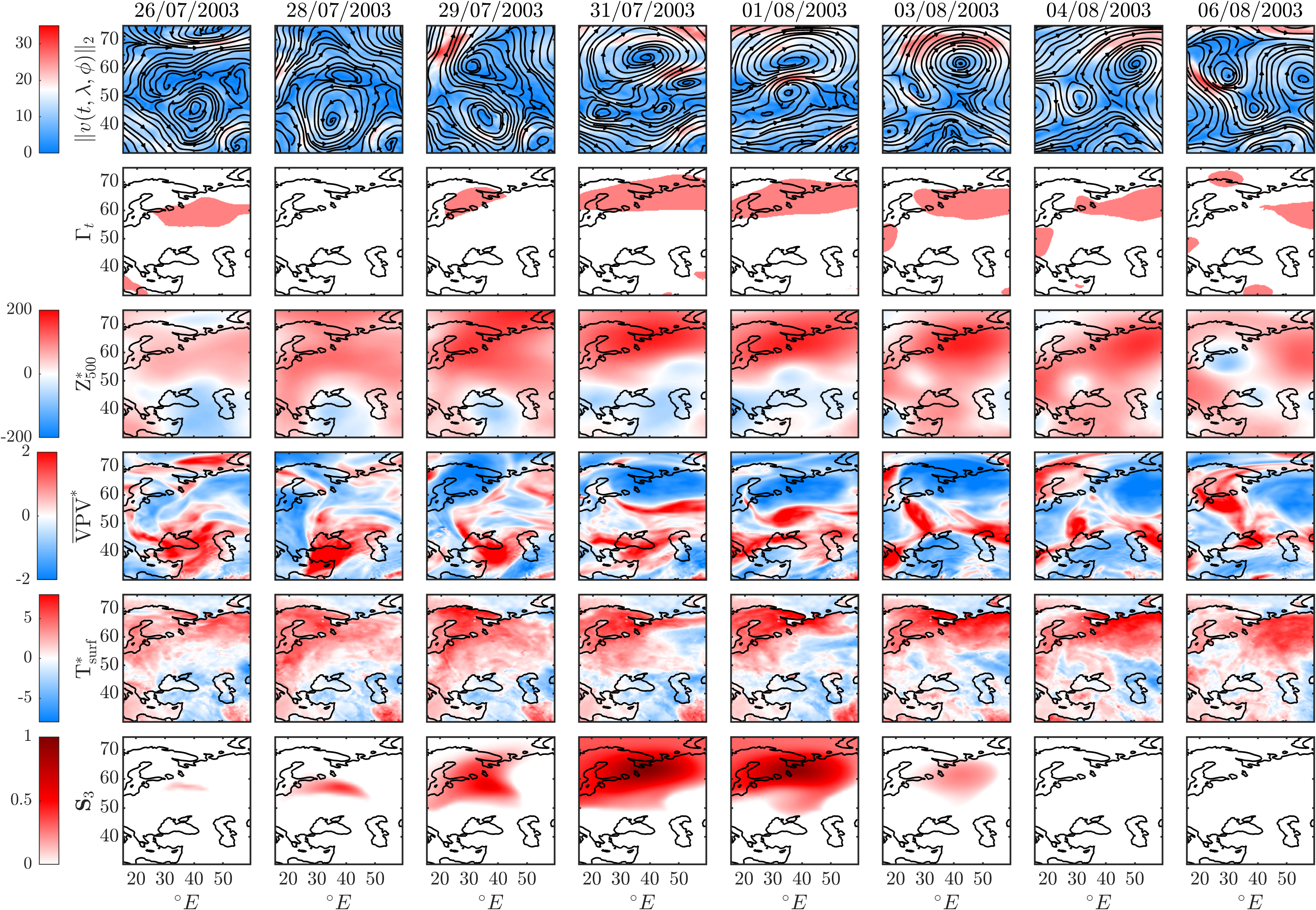}
    \caption{Several block diagnostic tools for the European atmospheric flow (eastern region) taken at specific dates over the date range 26/07/03-06/08/03. Included are (from top to bottom): the wind speed scalar field and corresponding streamlines; blocked regions $\Gamma_t$ identified using the $\mathrm{Z}_{500}$ gradient, the $\mathrm{Z}_{500}^{*}$ scalar field, the $\overline{\mathrm{VPV}}^{*}$ scalar field, \new{temperature anomalies taken at the Earth's surface $T_{\mathrm{surf}}^{*}$}, and the third SEBA vector $\mathbf{S}_3$ as shown in Figure \ref{fig:EuroBlock_East_InfGen}. The time of day is 00:00 (UTC) for all of the above dates. Blocks are identified through medium red coloured objects in the $\Gamma_t$ and $\mathrm{Z}_{500}^{*}$ fields; deep red objects in the $\mathbf{S}_3$ field; medium blue objects in the $\overline{\mathrm{VPV}}^{*}$ scalar field, and closed clockwise rotating streamline loops in the wind field.}
    \label{fig:EuroBlock_East_MultiFig}
\end{figure}

We now compare $\mathbf{S}_3$ to the geophysical block diagnostic methods listed in Section \ref{BlockIdentification}, as well as the wind streamline field \new{and surface temperature anomalies}. Figure \ref{fig:EuroBlock_East_MultiFig} shows time slices of $\mathbf{S}_3$ and all of the other diagnostic quantities taken at the same time instances as those used in Figure \ref{fig:EuroBlock_East_InfGen}.
All of the other block identification methods are able to identify the east block to some extent. This is evidenced by closed contours of clockwise oriented wind streamlines \cite{AnticycloneRotation}, blocked regions $\Gamma_t$ obtained from the $\mathrm{Z}_{500}$ gradient, and objects of strongly positive $\mathrm{Z}_{500}^{*}$ and negative $\overline{\mathrm{VPV}}^{*}$. \new{The block diagnosed is also reasonably correlated with strongly positive \emph{quasi-stationary} surface temperature anomalies (i.e.\ spatially correlated for those times where the temperature anomaly is relatively fixed in space).}
Regarding false positives, the $\mathrm{Z}_{500}$ gradient method shows a block on July 26--27 
%\aleks{In results not shown, yes. The blocked object disappears for one day only (July 28), it could be a separate (candidate) block that ends on July 27.}, 
that is not indicated by any of the other methods.
%None of the methods are able to identify any blocking activity before July 29 except for . %\gf{I see no evidence of this in the figure?} %\aleks{On July 26, a red object roughly positioned where the block is supposed to be is present in the $\Gamma_{t}$ field, on July 28 that object has vanished, and on July 29 another object has appeared slightly further West. None of the other diagnostic quantities show any evidence of a block being present on July 26, though I thought it would be prudent to mention that the gradient method finds a candidate blocked object on July 26, decidedly a false positive.}. 
After August 3, a block cannot be identified using $\mathbf{S}_3$.
In contrast, the results from each of the other methods give the impression that after August 3 the block moves Eastwards until August 6 where the block is roughly centred at $55^{\circ}$ E. This Eastward movement is not consistent with the quasi-stationary nature of blocking events \cite{Rex50I,Rex50II,Scherrer06} and we regard these temporal extensions of the block arising from the non-generator techniques as false positives. 
%From these results, the other block detection methods could be identifying a separate block present after August 4 and positioned further East as, in results not shown, the eigenvectors $\mathbf{F}_2^{*}$ and $\mathbf{F}_4^{*}$ as well as the SEBA vector $\mathbf{S}_1$ were found to detect another quasi-stationary object which forms on August 4 and behaves like a separate blocking event roughly longitudinally centred at $55^{\circ} \mathrm{E}$. 

We emphasise that the other blocking detection methods use instantaneous geophysical proxies to diagnose a block at individual times, and there is no explicit temporal connection between the diagnoses at distinct times. 
This disconnection in time not only ignores the important quasi-stationary aspect of blocks, but can make it difficult to identify separate blocking events and to piece together block identification across time, including when the block appears and disappears.
%In other words, a block can be identified at one time step using these indexing methods, but a completely different block can be found at a subsequent time step and, owing to that lack of temporal connection, there is no way to tell when one block has finished and another block has started simply by looking at these individual indication fields. 
In contrast, the inflated generator approach explicitly connects different time instances and enforces quasi-stationarity of the identified blocks, while the SEBA algorithm isolates individual blocks.
%we can clearly establish a more accurate picture of the temporal evolution of the precise blocking event we seek without interference from other quasi-stationary objects present in the flow; which is an outcome that cannot be achieved using the other diagnostic methods.
%The inflated generator method is not without its limitations. The biggest limitation is that we needed to have some idea of when the block is occurring and where it is positioned before identifying it using the generator. Eigenvectors of the generator were found to display several distinct quasi-stationary objects for this flow including the East block, a potential second block positioned further East, and a cut off low positioned below the block. Sometimes, combinations of these objects were identified through individual eigenvectors. We 
The inflated generator approach  requires some prior knowledge of the lifespan and spatial extent of a block to set up the spatial and temporal domain, but as one  can see from Figure \ref{fig:EuroBlock_East_MultiFig}, closed, clockwise rotating wind streamlines can provide this rough initial information.

\subsubsection{West Block (August 2003)} \label{EuroBlockWest}

Next we apply the inflated generator method to a block that occurred further west in Europe from 00:00 02/08/03 to 00:00 12/08/03 (UTC) (inclusive).
This significant blocking event was responsible for some of the worst extreme heat recorded in Western Europe (particularly in France) during the summer of 2003 \cite{Black04,Thomas21}. For this longer-duration block, we alter the spatial domain $M$ slightly to $\left[ 10^{\circ} \mathrm{W} , 40^{\circ} \mathrm{E} \right] \times \left[ 30^{\circ} \mathrm{N} , 75^{\circ} \mathrm{N} \right]$, and define a time interval covering the full duration of the block as well as five days immediately preceding and following this block, ranging from 00:00 28/07/03 to 00:00 17/08/03 (UTC) (inclusive), yielding $\tau = 20$ days. We maintain a temporal resolution of 6 hours and a spatial resolution of $1^{\circ}$ in both latitude and longitude, providing a spatial grid of $50 \times 45$ boxes.

\begin{comment}

\begin{figure}[h]
      \centering
      \includegraphics[width=0.9\textwidth]{EuroBlock_West_GeneratorEigVec_19Mar24.png}
      \caption{The leading non-trivial eigenvector $v_2$ and the second SEBA vector $S_2$ obtained from the dynamic generator for the European atmospheric flow (Western region) over the date range 28/07/03-17/08/03.}
      \label{fig:EuroBlock_West_DynGen}
\end{figure}

We begin by applying the dynamic generator method to this data in an attempt to identify the Western block. Figure \ref{fig:EuroBlock_West_DynGen} shows the scalar field for $v_2$, the first non--trivial eigenvector of the dynamic generator, as well as the scalar field for $S_2$, the second vector obtained after running the SEBA algorithm on the first two eigenvectors of the dynamic generator ($v_1$ (trivial) and $v_2$). The resultant SEBA vectors produced, $S_1$ and $S_2$, each have minimum values of 0 as desired. The coherent, quasi--stationary object identified by the dynamic generator is a small, semi--elliptic patch concentrated around North Western Africa and centred around $5^{\circ}$ longitude West and $35^{\circ}$ latitude North. This object does not extend into Europe, and does not cover the expected blocked region described earlier; which reinforces similar findings from the Double Gyre flow discussed in Section \ref{IdealisedDGFlow}.

\end{comment}

\begin{figure}[hbt!]
      \centering
      \includegraphics[width=0.5\textwidth]{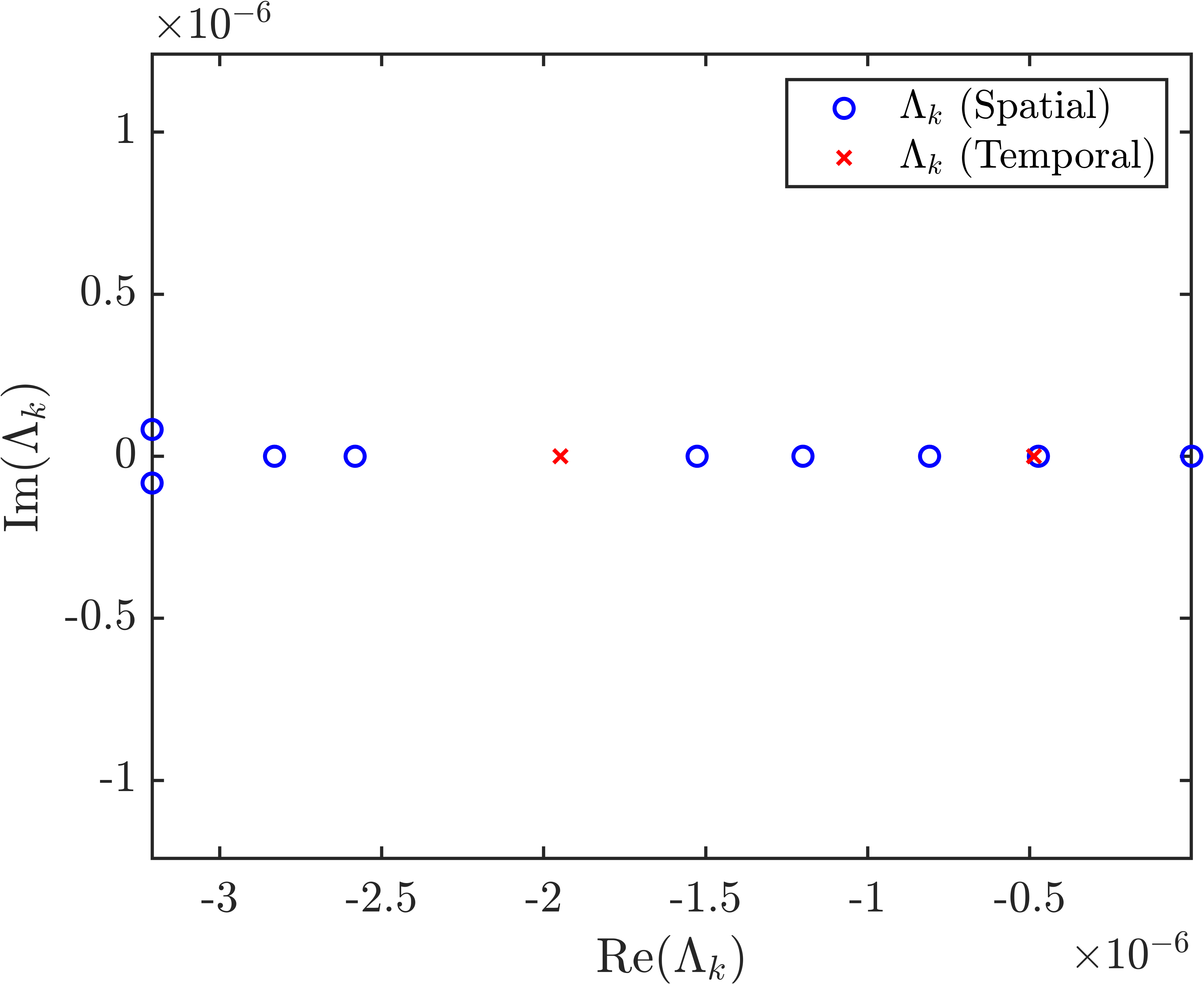}
      \caption{The leading eleven  eigenvalues $\Lambda_k$ (ordered according to their real parts in decreasing order) of the inflated generator calculated for the European atmospheric flow (western region) over the date range 28/07/03-17/08/03, plotted in the complex plane. Spatial and temporal eigenvalues are distinguished using separate plot markers and colours.}
      \label{fig:EuroBlock_West_Eigvals}
\end{figure}

\begin{figure}[hbt!]
      \centering
      \includegraphics[width=\textwidth]{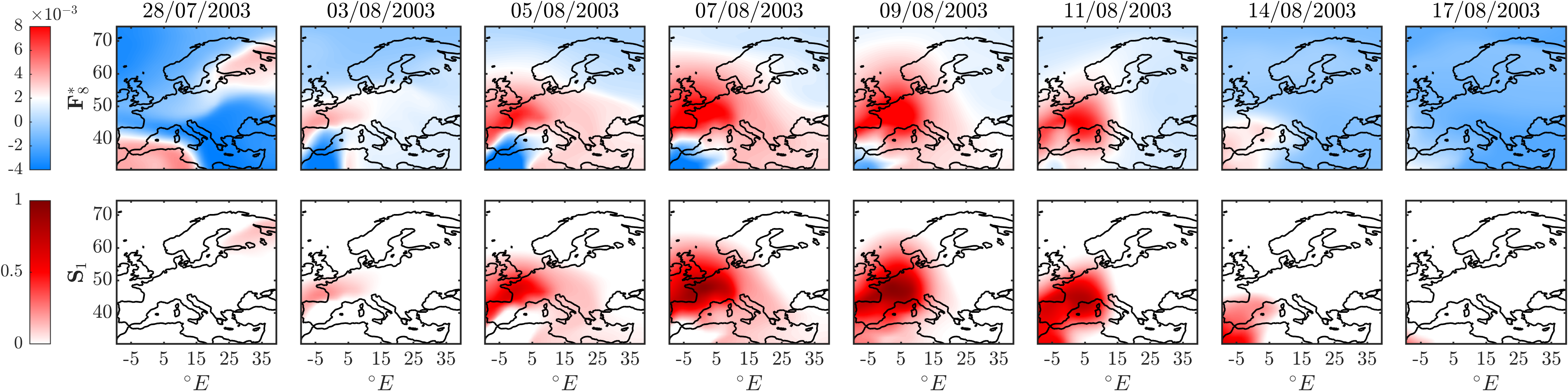}
      \caption{Temporal evolution of the spatial eigenvector $\mathbf{F}_8^{*}$ and the SEBA vector $\mathbf{S}_1$ of the inflated generator for the European atmospheric flow (western region) at select dates taken over the date range 28/07/03-17/08/03. The time of day is 00:00 (UTC) for all of the above dates.}
      \label{fig:EuroBlock_West_InfGen}
\end{figure}

We calculate $\epsilon \approx 327.5440 \ \mathrm{ms}^{-1/2}$ from a median box length $\bar{\ell} \approx 103618 \ \mathrm{m}$ and median speed $\bar{v} \approx 10.3539 \ \mathrm{ms^{-1}}$, and let $a = 0.0045$;  our heuristic yielded $0.0043$. The first eleven eigenvalues of the inflated generator spectrum are in Figure \ref{fig:EuroBlock_West_Eigvals}.
\new{The eigenvectors $\mathbf{F}_3^{*}$ and $\mathbf{F}_7^{*}$ are real valued and temporal, while $\mathbf{F}_{10}^{*}$ and $\mathbf{F}_{11}^{*}$ are complex valued;  the remaining eigenvectors are real valued and spatial.}
Similarly to the east block experiment, we inspect each non-trivial, real-valued spatial eigenvector until we find one that provides the best visualisation of the west block. %Initially, these eigenvectors do show the West block as part of a larger quasi-stationary object whose spread extends too far South West and takes its form before July 28. 
In this example we find that the eighth eigenvector $\mathbf{F}_8^{*}$ (see Figure \ref{fig:EuroBlock_West_InfGen} (upper row)) best isolates the west block as a separate quasi-stationary object covering
%\gf{covering means matching, subset or superset?} \aleks{I'm trying to say here that the block is positioned over France (as expected) over the anticipated time interval. I have altered the text.} \gf{Same question: What does "over" mean? It's like "cover", it's ambiguous.} 
the expected spatial and temporal extent of the block. 
The block is indicated by the deep red elliptical objects in the $\mathbf{F}_8^{*}$ fields, which form around August 4 and are positioned over  France between August 6-12, the period over which France's most extreme heat was recorded \cite{Black04,Thomas21}. 
The block is still present on August 12 and begins to disappear after this date before almost completely vanishing after August 14. 
To isolate this block from all other objects in the $\mathbf{F}_8^{*}$ field, we apply SEBA to all real-valued spatial eigenvectors from $\mathbf{F}_1^{*}$ up to and including $\mathbf{F}_8^{*}$.
%and inspect the six SEBA vectors generated to identify which of these best illustrates the West block. 
The first SEBA vector $\mathbf{S}_1$ isolates the west block and is plotted in Figure \ref{fig:EuroBlock_West_InfGen} (lower row).

\begin{figure}[hbt!]
    \centering
    \includegraphics[width=\textwidth]{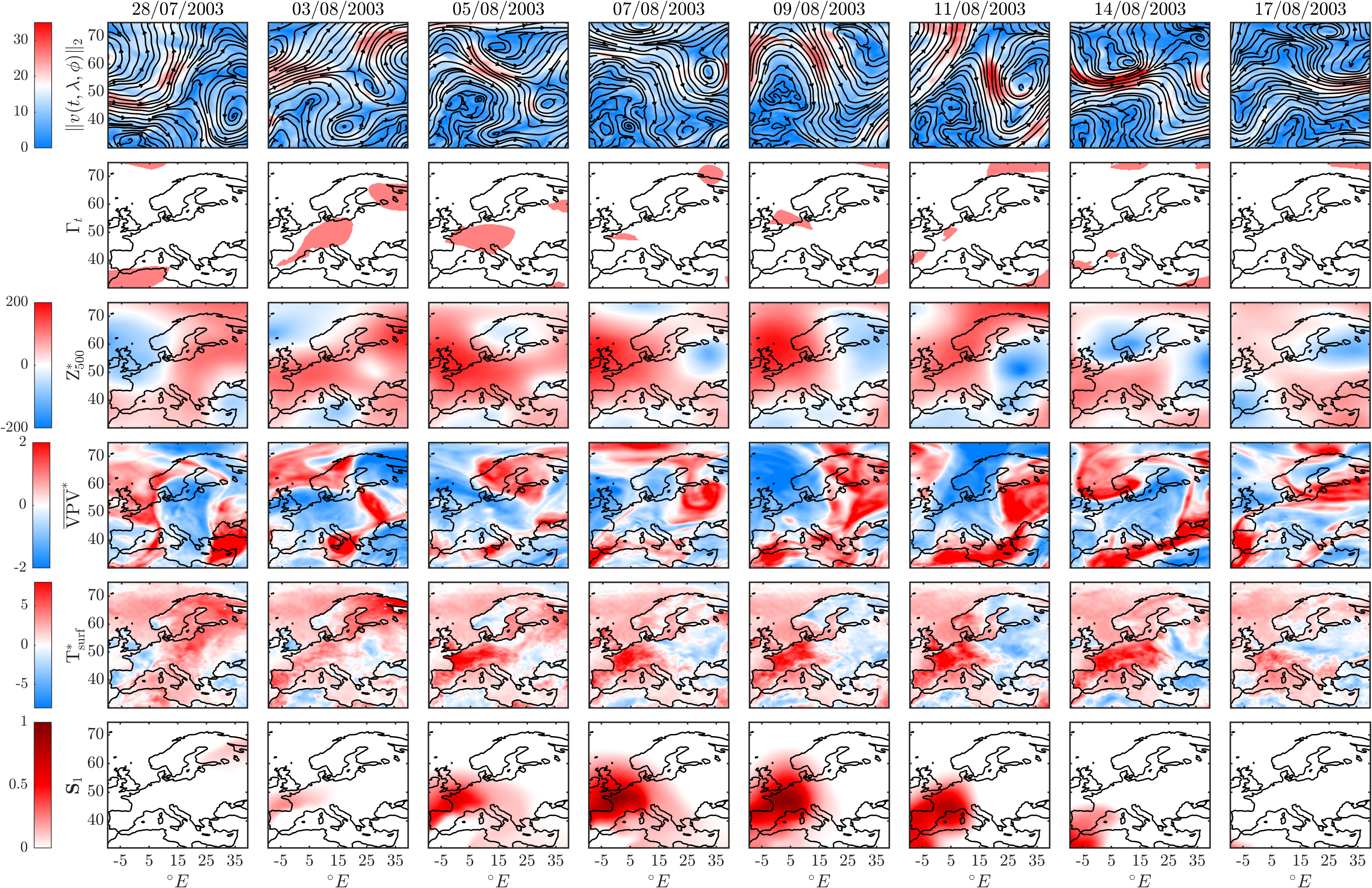}
    \caption{Several block diagnostic tools for the European atmospheric flow (western region) taken at select dates over the date range 28/07/03-17/08/03. Included are (from top to bottom): the wind speed scalar field and corresponding streamlines; blocked regions $\Gamma_t$ identified using the $\mathrm{Z}_{500}$ gradient, the $\mathrm{Z}_{500}^{*}$ scalar field, the $\overline{\mathrm{VPV}}^{*}$ scalar field, \new{temperature anomalies taken at the Earth's surface $T_{\mathrm{surf}}^{*}$},
     and the first SEBA vector $\mathbf{S}_1$ as shown in Figure \ref{fig:EuroBlock_West_InfGen}. The time of day is 00:00 (UTC) for all of the above dates. Blocks are identified through medium red coloured objects in the $\Gamma_t$ and $\mathrm{Z}_{500}^{*}$ fields; deep red objects in the $\mathbf{S}_1$ field; medium blue objects in the $\overline{\mathrm{VPV}}^{*}$ scalar field, and closed clockwise rotating streamline loops in the wind field.}
    \label{fig:EuroBlock_West_MultiFig}
\end{figure}

We now compare $\mathbf{S}_1$ to the other block diagnostics in Figure \ref{fig:EuroBlock_West_MultiFig}. 
The west block is considerably harder to identify using these geophysical proxies than the east block. 
%The closest any of these diagnostic tools get to matching the results identified through $\mathbf{S}_1$ are the wind streamlines, where a parabolic loop of clockwise oriented streamline contours \cite{AnticycloneRotation} centred over France has roughly formed by August 9. 
The $\mathrm{Z}_{500}^{*}$ and $\overline{\mathrm{VPV}}^{*}$ fields roughly identify a similar object to that identified using $\mathbf{S}_1$, though this object is centred over the United Kingdom and Ireland rather than France. The $\mathrm{Z}_{500}$ gradient struggles to identify the west block, initially identifying it on August 3, but only managing to identify random fragments of it from August 7 onwards. 
There appears to be no consistent agreement between the three geophysical block detection methods regarding the location, size and duration of the west block, and none were able to provide a clear identification of the west block. %or at the very least identify the block's correct position, shape and duration.

In contrast, the deep red elliptical object found in $\mathbf{S}_1$ strongly resembles the west block, with the largest values of $\mathbf{S}_1$ pinpointing the presence of a roughly elliptical structure over France during the expected period of the west block. 
\new{There is also a considerable correlation between this elliptical structure and quasi-stationary regions of strongly positive surface temperature anomalies.}
The paler red periphery of this object extends beyond the expected range of the block (particularly south west into Northwest Africa).
 %a key issue with the generator results for this block is the spatial extension of the identified quasi-stationary object further South West of France. 
 This extension of the west block was more strongly identifiable in all real-valued spatial eigenvectors preceding $\mathbf{F}_8^{*}$, suggesting that another blocking event was present over Northwest Africa between July 28 and August 3. Only the $\mathrm{Z}_{500}$ gradient on July 28 suggests the presence of this other block.
%Like with the East block, we needed to have some idea of the spatial and temporal extent of the West block in order to identify it using the inflated generator. The leading real-valued spatial eigenvector $\mathbf{F}_2^{*}$ again identified a combination of several quasi-stationary objects within the field, including the West block. We would need to use the wind streamlines to give an initial indication of the presence of a block, but we can then follow this up with the inflated generator method to develop a thorough picture of the onset, evolution and decay of this block. 

     \new{In the above two case studies, the sparse matrix $\mathbf{G}_a$ has of the order of $10^5$ rows and columns. Even at this size, computing a modest number of dominant eigenvalues and eigenvectors is still numerically efficient. In results not shown, at coarser levels of spatial and temporal resolution we were able to generate the same results as those shown for both blocking events with similar levels of detail.} %Additionally, our method shows satisfactory robustness against changes to the domain and boundary conditions used for each experiment. In other results not shown, we repeated the east and west block experiments with domains larger than those listed earlier for each block (maintaining Neumann boundary conditions), and using Dirichlet boundary conditions whilst maintaining the original domains. Despite some minor differences in the eigenvectors which contain each block, the SEBA vectors exhibiting each block remained stable compared to the results shown above.}

\begin{comment}
    This is owing to one of the $\zeta_{500}$ criteria (specifically $\gamma_S > 0$) not being satisfied over France at this particular period of time ($\gamma_S$ was found to be negative and of negligibly small value).
\end{comment}

\section{Conclusion}

We have developed a robust methodology for the identification of time-parameterised regions in phase space that are (i) almost-invariant or metastable at each time instant and (ii) vary slowly in space across time (are quasi-stationary).
Our approach also determined the timing of the onset and disappearance of such quasi-stationary families of almost-invariant sets from a single eigencomputation of an inflated generator.
In addition to introducing the concept of the inflated generator, we developed theory to underpin the almost-invariance and quasi-stationarity properties of the eigenvectors of the inflated generator, and proposed an efficient computational approach that requires only an estimate of the underlying time-varying velocity field (in particular, no trajectory integration is needed).
We illustrated our new approach by successfully identifying two atmospheric blocking events in the European 2003 summer, and compared our results with common blocking diagnosis techniques based on geophysical proxies.
We found that when compared to a ground truth, our inflated generator results were considerably more accurate than all geophysical proxies.
Advantages of our approach over existing geophysical methods for atmospheric block diagnosis include (i) the direct use of wind fields in assessing almost-invariance, rather than less direct proxies, and (ii) the explicit requirement that the diagnosed family of almost-invariant sets is quasi-stationary.
We hope that our spacetime methodology will aid a better understanding of the physics behind critical blocking phenomena through its complete lifecycle description of these events.

\section{Acknowledgements}

AB is supported by an ARC Discovery Project DP210100357, and GF is partially supported by an ARC Discovery Project DP210100357 and an ARC Laureate Fellowship FL230100088.
The authors acknowledge the use of UNSW's Katana compute cluster for some of the computations.

\section{Code}
The code used to produce the results featured in this article can be obtained from the Inflated-Generator repository on GitHub, available at: \url{https://github.com/gfroyland/Inflated-Generator}

%\section{(Appendix) Definition of Symbols/Abbreviations}
%\gf{not required}

%\nocite{*}
\bibliographystyle{abbrv}
\bibliography{citations}

\begin{thebibliography}{10}

\bibitem{AnomalyDefn}
{anomaly}.
\newblock \url{https://glossary.ametsoc.org/wiki/Anomaly}, 2024.
\newblock American Meteorological Society Glossary of Meteorology.

\bibitem{AnticycloneRotation}
{anticyclone}.
\newblock \url{https://glossary.ametsoc.org/wiki/Anticyclone}, 2024.
\newblock American Meteorological Society Glossary of Meteorology.

\bibitem{AnomalyComp}
{Climatologies and Standardized Anomalies}.
\newblock \url{https://iridl.ldeo.columbia.edu/dochelp/StatTutorial/Climatologies/index.html}, 2024.
\newblock Columbia University.

\bibitem{AFK24}
J.~Atnip, G.~Froyland, and P.~Koltai.
\newblock An inflated dynamic {L}aplacian to track the emergence and disappearance of semi-material coherent sets.
\newblock \url{https://arxiv.org/abs/2403.10360}, 2024.

\bibitem{Barnes12}
E.~A. Barnes, J.~Slingo, and T.~Woollings.
\newblock A methodology for the comparison of blocking climatologies across indices, models and climate scenarios.
\newblock {\em Climate Dynamics}, 38(11--12):2467--2481, 2012.

\bibitem{Barriopedro10}
D.~Barriopedro, R.~Garc{\'i}a-Herrera, and R.~Trigo.
\newblock Application of blocking diagnosis methods to {G}eneral {C}irculation {M}odels. {P}art {I}: a novel detection scheme.
\newblock {\em Climate Dynamics}, 35(1):1373--1391, 2010.

\bibitem{Black04}
E.~Black, M.~Blackburn, G.~Harrison, B.~Hoskins, and J.~Methven.
\newblock {Factors contributing to the summer 2003 European heatwave}.
\newblock {\em Weather (Royal Meteorological Society)}, 59(8):217--223, 2004.

\bibitem{Buehler11}
T.~Buehler, C.~C. Raible, and T.~F. Stocker.
\newblock {The relationship of winter season North Atlantic blocking frequencies to extreme cold or dry spells in the ERA-40}.
\newblock {\em Tellus A: Dynamic Meterology and Oceanography}, 63(2):212--222, 2011.

\bibitem{Cai20}
W.~Cai, X.~Xu, X.~Cheng, F.~Wei, X.~Qiu, and W.~Zhu.
\newblock Impact of ``blocking" structure in the troposphere on the wintertime persistent heavy air pollution in {N}orthern {C}hina.
\newblock {\em Science of the Total Environment}, 741(1):140325, 2020.

\bibitem{Carrera04}
M.~Carrera, R.~Higgins, and V.~Kousky.
\newblock {Downstream Weather Impacts Associated with Atmospheric Blocking over the Northeast Pacific}.
\newblock {\em Journal of Climate}, 17(24):4823--4839, 2004.

\bibitem{DJ99}
M.~Dellnitz and O.~Junge.
\newblock {On the Approximation of Complicated Dynamical Behavior}.
\newblock {\em SIAM Journal on Numerical Analysis}, 36(2):491--515, 1999.

\bibitem{Deuflhard05}
P.~Deuflhard and M.~Weber.
\newblock {Robust Perron cluster analysis in conformation dynamics}.
\newblock {\em Linear Algebra and its Applications}, 398(1):161--184, 2005.

\bibitem{Dole83}
R.~M. Dole and N.~D. Gordon.
\newblock {Persistent Anomalies of the Extratropical Northern Hemisphere Wintertime Circulation: Geographical Distribution and Regional Persistence Characteristics}.
\newblock {\em Monthly Weather Review}, 111(8):1567--1586, 1983.

\bibitem{DunnSigouin13}
E.~Dunn-Sigouin and S.-W. Son.
\newblock Northern {H}emisphere blocking frequency and duration in the {CMIP}5 models.
\newblock {\em Journal of Geophysical Research: Atmospheres}, 118(3):1179--1188, 2013.

\bibitem{Ehstand21}
N.~Ehstand, R.~V. Donner, C.~L\'{o}pez, and E.~Hern\'{a}ndez-Garc\'{i}a.
\newblock Characteristic signatures of {N}orthern {H}emisphere blocking events in a {L}agrangian flow network representation of the atmospheric circulation.
\newblock {\em Chaos: An Interdisciplinary Journal of Nonlinear Science}, 31(9):093128, 2021.

\bibitem{Faranda15}
D.~Faranda, G.~Masato, N.~Moloney, Y.~Sato, F.~Daviaud, B.~Dubrulle, and P.~Yiou.
\newblock {The switching between zonal and blocked mid-latitude atmospheric circulation: a dynamical system perspective}.
\newblock {\em Climate Dynamics}, 47(5--6):1587--1599, 2015.

\bibitem{F05}
G.~Froyland.
\newblock Statistically optimal almost--invariant sets.
\newblock {\em Physica D: Nonlinear Phenomena}, 200(3--4):205--219, 2005.

\bibitem{F13}
G.~Froyland.
\newblock An analytic framework for identifying finite-time coherent sets in time-dependent dynamical systems.
\newblock {\em Physica D: Nonlinear Phenomena}, 250(1):1--19, 2013.

\bibitem{FD03}
G.~Froyland and M.~Dellnitz.
\newblock Detecting and locating near-optimal almost-invariant sets and cycles.
\newblock {\em SIAM Journal on Scientific Computing}, 24(6):1839--1863, 2003.

\bibitem{FJK13}
G.~Froyland, O.~Junge, and P.~Koltai.
\newblock Estimating long-term behavior of flows without trajectory integration: The infinitesimal generator approach.
\newblock {\em SIAM Journal on Numerical Analysis}, 51(1):223--247, 2013.

\bibitem{FK23}
G.~Froyland and P.~Koltai.
\newblock Detecting the birth and death of finite-time coherent sets.
\newblock {\em Communications on Pure and Applied Mathematics}, 76(12):3642--3684, 2023.

\bibitem{FKS20}
G.~Froyland, P.~Koltai, and M.~Stahn.
\newblock {Computation and Optimal Perturbation of Finite-Time Coherent Sets for Aperiodic Flows Without Trajectory Integration}.
\newblock {\em SIAM Journal on Applied Dynamical Systems}, 19(3):1659--1700, 2020.

\bibitem{FLS10}
G.~Froyland, S.~Lloyd, and N.~Santitissadeekorn.
\newblock {Coherent sets for nonautonomous dynamical systems}.
\newblock {\em Physica D: Nonlinear Phenomena}, 239(16):1527--1541, 2010.

\bibitem{FP09}
G.~Froyland and K.~Padberg.
\newblock {Almost--invariant sets and invariant manifolds — Connecting probabilistic and geometric descriptions of coherent structures in flows}.
\newblock {\em Physica D: Nonlinear Phenomena}, 238(16):1507--1523, 2009.

\bibitem{FRS19}
G.~Froyland, C.~P. Rock, and K.~Sakellariou.
\newblock {Sparse eigenbasis approximation: Multiple feature extraction across spatiotemporal scales with application to coherent set identification}.
\newblock {\em Communications in Nonlinear Science and Numerical Simulation}, 77(1):81--107, 2019.

\bibitem{FSM10}
G.~Froyland, N.~Santitissadeekorn, and A.~Monahan.
\newblock {Transport in time-dependent dynamical systems: Finite--time coherent sets}.
\newblock {\em Chaos: An Interdisciplinary Journal of Nonlinear Science}, 20(4):043116, 2010.

\bibitem{FSvS14}
G.~Froyland, R.~M. Stuart, and E.~van Sebille.
\newblock {How well-connected is the surface of the global ocean?}
\newblock {\em Chaos: An Interdisciplinary Journal of Nonlinear Science}, 24(3):033126, 2014.

\bibitem{Ghil18}
M.~Ghil, A.~Groth, D.~Kondrashov, and A.~W. Robinson.
\newblock {Extratropical sub-seasonal–to–seasonal oscillations and multiple regimes: The dynamical systems view}.
\newblock In A.~W. Robinson and F.~Vitart, editors, {\em {Sub-Seasonal to Seasonal Prediction: The Gap Between Weather and Climate Forecasting}}, chapter~6, pages 119--142. Elsevier, Amsterdam, 2018.

\bibitem{HBV13}
G.~Haller and F.~Beron-Vera.
\newblock {Coherent Lagrangian vortices: the black holes of turbulence}.
\newblock {\em Journal of Fluid Mechanics}, 731(1):R4, 2013.

\bibitem{HP98}
G.~Haller and A.~Poje.
\newblock {Finite time transport in aperiodic flows}.
\newblock {\em Physica D: Nonlinear Phenomena}, 119(3--4):352--380, 1998.

\bibitem{Haller00}
G.~Haller and G.~Yuan.
\newblock {Lagrangian coherent structures and mixing in two–-dimensional turbulence}.
\newblock {\em Physica D: Nonlinear Phenomena}, 147(3--4):352--370, 2000.

\bibitem{Hauser23}
S.~Hauser, F.~Teubler, M.~Riemer, P.~Knippertz, and C.~M. Grams.
\newblock {Towards a holistic understanding of blocked regime dynamics through a combination of complementary diagnostic perspectives}.
\newblock {\em Weather and Climate Dynamics}, 4(2):399--425, 2023.

\bibitem{ERA5}
H.~Hersbach, B.~Bell, P.~Berrisford, G.~Biavati, A.~Hor\'{a}nyi, J.~Muñoz-Sabater, J.~Nicolas, C.~Peubey, R.~Radu, I.~Rozum, D.~Schepers, A.~Simmons, C.~Soci, D.~Dee, and J.-N. Th\'{e}paut.
\newblock {ERA5 hourly data on pressure levels from 1940 to present, Copernicus Climate Change Service (C3S) Climate Data Store (CDS)}.
\newblock \url{https://doi.org/10.24381/cds.bd0915c6}, 2023.

\bibitem{ERA5_singlelevels}
H.~Hersbach, B.~Bell, P.~Berrisford, G.~Biavati, A.~Hor\'{a}nyi, J.~Muñoz-Sabater, J.~Nicolas, C.~Peubey, R.~Radu, I.~Rozum, D.~Schepers, A.~Simmons, C.~Soci, D.~Dee, and J.-N. Th\'{e}paut.
\newblock {ERA5 hourly data on single levels from 1940 to present, Copernicus Climate Change Service (C3S) Climate Data Store (CDS)}.
\newblock \url{https://doi.org/10.24381/cds.adbb2d47}, 2023.

\bibitem{Hirt18}
M.~Hirt, L.~Schielicke, A.~M\"{u}ller, and P.~N\'{e}vir.
\newblock {Statistics and dynamics of blockings with a point vortex model}.
\newblock {\em Tellus A: Dynamic Meteorology and Oceanography}, 70(1):1458565, 2018.

\bibitem{quadgk}
S.~G. Johnson.
\newblock {QuadGK.jl}: {G}auss--{K}ronrod integration in {J}ulia: v2.9.4.
\newblock \url{https://github.com/JuliaMath/QuadGK.jl}, 2013.

\bibitem{Kautz22}
L.-A. Kautz, O.~Martius, S.~Pfahl, J.~G. Pinto, A.~M. Ramos, P.~M. Sousa, and T.~Woollings.
\newblock {Atmospheric blocking and weather extremes over the Euro-Atlantic sector – a review}.
\newblock {\em Weather and Climate Dynamics}, 3(1):305--336, 2022.

\bibitem{Interpolations}
M.~Kittisopikul, T.~E. Holy, and T.~Aschan.
\newblock {JuliaMath/Interpolations.jl: v0.14.7}.
\newblock \url{https://github.com/JuliaMath/Interpolations.jl}, 2022.

\bibitem{K11}
P.~Koltai.
\newblock {A stochastic approach for computing the domain of attraction without trajectory simulation}.
\newblock {\em American Institute of Mathematical Sciences Conference Publications}, 2011(Special):854--863, 2011.

\bibitem{Lejenas83}
H.~Lejen\"{a}s and H.~{\O}kland.
\newblock Characteristics of northern hemisphere blocking as determined from a long time series of observational data.
\newblock {\em Tellus A: Dynamic Meterology and Oceanography}, 35(5):350--362, 1983.

\bibitem{Liu94}
Q.~Liu.
\newblock On the definition and persistence of blocking.
\newblock {\em Tellus A: Dynamic Meterology and Oceanography}, 46(3):286--298, 1994.

\bibitem{MB14}
T.~Ma and E.~M. Bollt.
\newblock {Differential geometry perspective of shape coherence and curvature evolution by finite-time nonhyperbolic splitting}.
\newblock {\em SIAM Journal on Applied Dynamical Systems}, 13(3):1106--1136, 2014.

\bibitem{mcowen}
R.~C. McOwen.
\newblock {\em {Partial Differential Equations: Methods and Applications}}.
\newblock Pearson, London, first edition, 1995.

\bibitem{Nakamura18}
N.~Nakamura and C.~S.~Y. Huang.
\newblock {Atmospheric blocking as a traffic jam in the jet stream}.
\newblock {\em Science}, 361(6397):42--47, 2018.

\bibitem{Park11}
T.-W. Park, C.-H. Ho, and S.~Yang.
\newblock {Relationship between the Arctic Oscillation and Cold Surges over East Asia}.
\newblock {\em Journal of Climate}, 24(1):68--83, 2011.

\bibitem{Pelly03}
J.~L. Pelly and B.~J. Hoskins.
\newblock {A New Perspective on Blocking}.
\newblock {\em Journal of the Atmospheric Sciences}, 60(5):743--755, 2003.

\bibitem{PY93}
R.~T. Pierrehumbert and H.~Yang.
\newblock {Global chaotic mixing on isentropic surfaces}.
\newblock {\em Journal of the Atmospheric Sciences}, 50(15):2462--2480, 1993.

\bibitem{Pinheiro19}
M.~C. Pinheiro, P.~A. Ullrich, and R.~Grotjahn.
\newblock Atmospheric blocking and intercomparison of objective detection methods: flow field characteristics.
\newblock {\em Climate Dynamics}, 53(7):4189--4216, 2019.

\bibitem{Rex50I}
D.~F. Rex.
\newblock Blocking {A}ction in the {M}iddle {T}roposphere and its {E}ffect upon {R}egional {C}limate {P}art {I}: {A}n {A}erological {S}tudy of {B}locking {A}ction.
\newblock {\em Tellus A: Dynamic Meterology and Oceanography}, 2(3):196--211, 1950.

\bibitem{Rex50II}
D.~F. Rex.
\newblock Blocking {A}ction in the {M}iddle {T}roposphere and its {E}ffect upon {R}egional {C}limate {P}art {II}: {T}he {C}limatology of {B}locking {A}ction.
\newblock {\em Tellus A: Dynamic Meterology and Oceanography}, 2(4):237--324, 1950.

\bibitem{RKLW90}
V.~Rom-Kedar, A.~Leonard, and S.~Wiggins.
\newblock {An analytical study of transport, mixing and chaos in an unsteady vortical flow}.
\newblock {\em Journal of Fluid Mechanics}, 214(1):347--394, 1990.

\bibitem{Rowley09}
C.~W. Rowley, I.~Mezi\'{c}, S.~Bagheri, P.~Schlatter, and D.~S. Henningson.
\newblock {Spectral analysis of nonlinear flows}.
\newblock {\em Journal of Fluid Mechanics}, 641(1):115--127, 2009.

\bibitem{Sausen95}
R.~Sausen, W.~K\"{o}nig, and F.~Sielmann.
\newblock Analysis of blocking events from observations and {ECHAM} model simulations.
\newblock {\em Tellus A: Dynamic Meterology and Oceanography}, 47(4):421--438, 1995.

\bibitem{Scherrer06}
S.~C. Scherrer, M.~Croci-Maspoli, C.~Schwierz, and C.~Appenzeller.
\newblock {Two-dimensional indices of atmospheric blocking and their statistical relationship with winter climate patterns in the Euro-Atlantic region}.
\newblock {\em International Journal of Climatology}, 26(2):233--249, 2006.

\bibitem{Schoeller24}
H.~Schoeller, R.~Chemnitz, P.~Koltai, M.~Engel, and S.~Pfahl.
\newblock {Assessing Lagrangian Coherence in Atmospheric Blocking}.
\newblock \url{https://doi.org/10.5194/egusphere-2024-2173}, 2024.

\bibitem{Schubert16}
S.~Schubert and V.~Lucarini.
\newblock {Dynamical analysis of blocking events: spatial and temporal fluctuations of covariant Lyapunov vectors}.
\newblock {\em Quarterly Journal of the Royal Meteorological Society}, 142(698):2143--2158, 2016.

\bibitem{SchuetteHabilitation}
C.~Sch\"{u}tte.
\newblock {\em {Conformational Dynamics: Modelling, Theory, Algorithm, and Application to Biomolecules}}.
\newblock Habilitation thesis, Zuse Institute, Freie Universit\"{a}t Berlin, July 1999.

\bibitem{Schwierz04}
C.~Schwierz, M.~Croci-Maspoli, and H.~C. Davies.
\newblock Perspicacious indicators of atmospheric blocking.
\newblock {\em Geophysical Research Letters}, 31(6):L06125, 2004.

\bibitem{SerGiacomi15}
E.~Ser-Giacomi, R.~Vasile, I.~Recuerda, E.~Hern\'{a}ndez-Garc\'{i}a, and C.~L\'{o}pez.
\newblock {Dominant transport pathways in an atmospheric blocking event}.
\newblock {\em Chaos: An Interdisciplinary Journal of Nonlinear Science}, 25(8):087413, 2015.

\bibitem{Shaddenetal05}
S.~C. Shadden, F.~Lekien, and J.~E. Marsden.
\newblock {Definition and properties of Lagrangian coherent structures from finite--time Lyapunov exponents in two--dimensional aperiodic flows}.
\newblock {\em Physica D: Nonlinear Phenomena}, 212(3--4):271--304, 2005.

\bibitem{Shutts83}
G.~J. Shutts.
\newblock {The propagation of eddies in diffluent jetstreams: Eddy vorticity forcing of `blocking' flow fields}.
\newblock {\em Quarterly Journal of the Royal Meteorological Society}, 109(462):737--761, 1983.

\bibitem{Steinfeld19}
D.~Steinfeld and S.~Pfahl.
\newblock {The role of latent heating in atmospheric blocking dynamics: a global climatology}.
\newblock {\em Climate Dynamics}, 53(9--10):6159--6180, 2019.

\bibitem{Tantet15}
A.~Tantet, F.~R. van~der Burgt, and H.~A. Dijkstra.
\newblock {An early warning indicator for atmospheric blocking events using transfer operators}.
\newblock {\em Chaos: An Interdisciplinary Journal of Nonlinear Science}, 25(3):036406, 2015.

\bibitem{Thomas21}
C.~Thomas, A.~Voulgarakis, G.~Lim, J.~Haigh, and P.~Nowack.
\newblock An unsupervised learning approach to identifying blocking events: the case of {E}uropean summer.
\newblock {\em Weather and Climate Dynamics}, 2(3):581--608, 2021.

\bibitem{Tibaldi90}
S.~Tibaldi and F.~Molteni.
\newblock On the operational predictability of blocking.
\newblock {\em Tellus A: Dynamic Meterology and Oceanography}, 42(3):343--365, 1990.

\bibitem{Tibaldi95}
S.~Tibaldi, P.~Ruti, E.~Tosi, and M.~Maruca.
\newblock Operational predictability of winter blocking at {ECMWF}: an update.
\newblock {\em Annales Geophysicae}, 13(1):305--317, 1995.

\bibitem{VonLindheim21}
J.~von Lindheim, A.~Harikrishnan, T.~D\"{o}rffel, R.~Klein, P.~Koltai, N.~Mikula, A.~M\"{u}ller, P.~N\'{e}vir, G.~Pacey, R.~Polzin, and N.~Vercauteren.
\newblock Definition, detection, and tracking of persistent structures in atmospheric flows.
\newblock {\em Journal of Mathematical Physics}, pages 1--64, 2021.

\bibitem{Xavier24}
A.~K. Xavier, J.~Demaeyer, and S.~Vannitsem.
\newblock {Variability and predictability of a reduced-order land–atmosphere coupled model}.
\newblock {\em Earth System Dynamics}, 15(4):893--912, 2024.

\bibitem{Yun19}
S.~Yun and C.~Yoo.
\newblock {The Effects of Spring and Winter Blocking on PM10 Concentration in Korea}.
\newblock {\em Atmosphere}, 10(7):410, 2019.

\end{thebibliography}

\end{document}